\documentclass{amsart}

  \usepackage[all]{xy}

  \usepackage{epsf,epsfig,amsfonts,graphicx,color}
  
  \usepackage{caption}
\usepackage{subcaption}

\usepackage{amsmath,amssymb}
  \usepackage{url}
   
 \def\R{\mathbb R}

\def\eps{\epsilon}

 \def\on{ orthonormal }

\def\Ri{ Riemannian }
\def\sR{ subRiemannian }

\def\g{\mathfrak{g}}

\newcommand{\beq}{\begin{equation}}
\newcommand{\eeq}{\end{equation}}

\newcommand{\dd}[2] { {{\partial #1}   \over {\partial #2}} }

\newtheorem{theorem}{Theorem}

\newtheorem{proposition}{Proposition}[section]

\newtheorem{lemma}{Lemma}[section]
\newtheorem{definition}{Definition}[section]
\newtheorem{remark}{Remark}[section]

\begin{document}

\newtheorem*{backgroundtheorem}{Background Theorem}


\title[Geodesics in Jet Space.]{Geodesics  in  Jet Space.}  
\author[A.\ Bravo-Doddoli]{Alejandro\ Bravo-Doddoli} \author[R.\ Montgomery]{Richard Montgomery}
\address{Alejandro Bravo-Doddoli: Dept. of Mathematics, UCSC,
1156 High Street, Santa Cruz, CA 95064}
\email{Abravodo@ucsc.edu}
\address{R. Montgomery: Dept. of Mathematics, UCSC,
1156 High Street, Santa Cruz, CA, 95064}
\email{rmont@ucsc.edu}
\begin{abstract}
The space $J^k$ of $k$-jets of a real function of one real variable $x$  admits the  structure of Carnot group type.   As such, $J^k$ admits
a   submetry (\sR submersion) onto the Euclidean plane. Horizontal lifts of  Euclidean lines (which
are the left-translates of  horizontal one-parameter subgroups) are thus  globally minimizing geodesics on $J^k$.  
All $J^k$-geodesics, minimizing or not, are constructed from degree $k$ polynomials in $x$ according to \cite{Monroy1},\cite{Monroy2},\cite{Monroy3},
reviewed here.  
The constant polynomials correspond to the horizontal lifts of lines. Which other polynomials yield globally minimizers and what
do these minimizers look like? We give a partial answer.    Our methods include constructing an intermediate three-dimensional ``magnetic'' \sR space  lying between the jet space and the plane,  solving a  Hamilton-Jacobi (eikonal) equations on this space,  and analyzing  period  asymptotics associated to  period degenerations arising from two-parameter families of these polynomials.  Along the way, we   conjecture   the  independence  of the  cut time of any geodesic on jet space  from the starting location on that geodesic.  
\end{abstract} 
\thanks{ }

\maketitle

\section{Introduction: Motivation, results, acknowledgement}

It is a basic and important fact that lines in Euclidean space are globally minimizing geodesics. 
Not only are   lines   geodesics, but no matter how far out we  travel   along a line away from a  point on the line,
the corresponding line segment continues to minimize the   distance between its end points.     Contrast this with the case of
geodesics on a cylinder, where most geodesics eventually fail to be minimizing.
     In the context of Carnot groups we can write down 
geodesic equations which describe  most  geodesics.  (They miss the ``abnormal'' or ``singular geodesics''.
See \cite{tour}.)  The horizontal lines -- the  left translates of horizontal one-parameter subgroups -   are 
 globally minimizing geodesics. 
In the first non-trivial case,  the Heisenberg group, the horizontal lines exhaust the set of  globally minimizing geodesics.
What happens for other Carnot groups?   
Are there any other globally minimizing
geodesics besides the horizontal lines ?  

The spaces $J^k = J^k (\R, \R)$ of $k$-jets of a real function  of a single real variable  forms
a   family of $k+2$-dimensional Carnot groups. (See \cite{CarnotJets}.)  $J^k$ is  the unique Carnot group of its dimension 
Goursat type: its  Lie bracket growth vector  is   $(2,3, 4, \ldots, k+2)$.  
$J^1$ is the   well-known Heisenberg group and, as we just saw, has no global minimizers
beyond the horizontal lines.  $J^2$ is the   Engel group \cite{tour}  
and has exactly one new global minimizer up to translation and scaling,   this geodesic being the horizontal lift of   the  ``Euler soliton'' whose global minimality is established in   \cite{Ardentov1,Ardentov2}.  See the middle panel of  figure \ref{fig:kink}.

 Anzaldo-Meneses and Monroy-Per\'ez  \cite{Monroy1,Monroy2,Monroy3}  showed that the \sR  geodesic flow on  $J^k$ is completely integrable.  In doing so they    parameterized the space of all geodesics (modulo Carnot translations)
 by an open subset of the  space of real polynomials $F(x)$ of degree $k$ (modulo  translation $F(x) \mapsto F(x-x_0)$).
 
Re-iterating, a  geodesic is called    globally minimizing if
each of its  compact subarcs  realizes the distance between its endpoints.
Our goal in this paper  is to select out those  degree $k$ polynomials which  yield  global minimizers on $J^k$, $k > 2$.

We  partially succeed. 
Theorem \ref{theorem-1} below excludes most polynomials from yielding global minimizers.  Theorem \ref{main} establishes the existence of 
 a previously unknown 8-dimensional family of global minimizers and characterizes  them  in terms of their   polynomials. 
 These two theorems are described in the next section, section \ref{sec: setup}.    The question of finding an exact characterization of the global minimizers in terms of their polynomials
remains open.     

Our methods  are three-fold.  First, in section \ref{sec: magnetic}, for each choice of polynomial $F(x)$ we construct an intermediate 3-dimensional \sR ``magnetic space'' denoted $\R^3 _F$ which 
lies between $J^k$ and the Euclidean plane and we reduce most of our work to analysis on this space.
Second, in section \ref{sec-H-J} we apply a  Hamilton-Jacobi method (also known as the method of calibrations)  to insure  that our candidate globally minimizing geodesics actually globally  minimize within a large
open slab-like domain   which contains them.  Finally we are reduced to a detailed
analysis of all the geodesics in the magnetic space which leave the slab-like domain to finish off the proof. 
In this last (exhausting) step which takes up section \ref{sec:proof-main-the}  we show that none of these competitor geodesics are  simultaneously   shorter and match endpoint conditions with our candidate geodesics.

\subsection{Acknowledgement}

It is an honor and a sadness to put forth this article  
in a Journal issue dedicated to Alexey Borisov.  In addition to stating our gratitude to Borisov and condolences
to his family and friends we would like to    thank Andrei Ardentov,  Gil Bor, Eero Hakavouri, Enrico Le Donne,  Hector Sanchez-Morgado, 
Felipe Monroy-Perez,   and Yuri Sachkov for   e-mail conversations regarding  the course of this work. We would also like to thank the
three anonymous reviewers for their diligent work and useful suggestions.
 This paper was developed with the support of the scholarship (CVU 619610) from  "Consejo de Ciencia y Tecnologia"  (CONACYT).

\vskip .3cm 

\section{ Set-up. Background. Theorems.   Overview.}\label{sec: setup}
 \subsection{Set-up and Carnot group structure.} 

We say that smooth real-valued function $f(x)$ and $g(x)$ are equivalent up to order $k$  at $x_0$
if $f(x) - g(x) = O(|x-x_0|^{k+1})$ holds.  Being equivalent to order $k$ is an equivalence relation on the space of 
germs of  smooth functions at $x_0$ 
and an equivalence class is called a k-jet at $x_0$.  The $k$-jet of a function  $f$ at $x_0$  can be identified with its $k$th order
Taylor expansion of $f$ at $x_0$ and,  as such,  is determined by the list of its first $k$ derivatives at $x_0$: 
$$ u_0  = f(x_0), \qquad u_j = d^j f / dx ^j (x_0):= f^{(j)} (x_0), j = 1, \ldots , k.$$ By letting the base point and function vary
we sweep out  the $k$-jet space $J^k$, a $k+2$-dimensional manifold with global coordinates $x$ and these $u_j$'s.  

If we fix the function $f$ and let the independent variable  $x$ vary, we get a curve   $j^k f: \R \to J^k$  called  the $k$-jet of $f$,  sending
$x \in \R$ to  the $k$-jet of $f$ at $x$. In coordinates 
$$(j^k f ) (x) = (x,  u_k (x), u_{k-1} (x), \ldots , u_1 (x), u_0 (x)) ;\;\; u_i (x)  = f^{(i)} (x)).$$ 
The  $k$-jet curve itself    is  everywhere tangent to the  rank two  distribution $D \subset TJ^k$ 
which is  globally framed by the two vector fields
 \begin{equation}
 \label{frame}
X_1 = \frac{\partial}{\partial x}  + \sum_{i=1}^{k} u_i \frac{\partial}{\partial u_{i-1}} \;\; \text{and} \;\; X_2 = \frac{\partial}{\partial u_k}. 
\end{equation}
A  \sR structure on $J^k$  is defined by declaring these two vector  fields to be \on.
In coordinates the \sR metric is defined by restricting $ds^2 = dx^2+du_k^2$ to $D$. 
Now  
$$\frac{d}{dx} j^k f (x) = X_1 + f^{(k+1)} (x) X_2$$
so that the  \sR length $\ell$   of the  curve $x \mapsto j^k f (x)$, restricted to a  
 finite interval $a \le x \le b$  is  
 $$\ell(j^k f |_{[a,b]})= \int_a ^b \sqrt{ 1 + (f^{(k+1)}(x))^2  }dx.$$

The map $\pi: J^k \to \R^2$ defined by
$$\pi(x, u_k, u_{k-1}, \ldots , u_0) = (x, u_k)$$
defines a   \sR submersion (or submetry) onto the Euclidean plane.  In other words,
its restriction to each two-plane $D$ is an isometry onto the
Euclidean plane with Euclidean metric $dx^2+du_k^2$.
This projection has an `inverse map', the horizontal lift, on the level of curves.
To understand the lift, rewrite $D$ as a Pfaffian system:
\begin{equation}
 \label{Pfaff}
 \begin{aligned}
 du_{i-1} - u_idx& = & 0\;\;\; \text{with}\; 1 \leq i \leq k. 
 \end{aligned}
  \end{equation}
  For example, the   last  of these equations, the one for $i = k$,  reads $du_{k-1} = u_k dx$.
 Given a smooth curve $c(t) = (x(t), u_k (t))$ in the plane 
 we associate to it the following   {\it horizontal lift equations}  \begin{equation}
\begin{aligned}
\label{eq:horiz} 
\dot u_{i} (t) & = u_{i+1} (x) \dot x(t),  \qquad  i =0, \ldots , k-1
\end{aligned}
\end{equation}
which simply say that the curve $\gamma(t) = (x(t), u_k (t), u_{k-1} (t), \ldots , u_0 (t))$
is horizontal and projects onto $c(t)$.   We call these curves $\gamma$ the
horizontal lifts of our plane curve.   The length of $\gamma$ and of $c$ over any compact
time interval are equal.   The horizontal lift $\gamma(t)$
is uniquely specified by the choice of initial condition, say $\gamma(0)$,  corresponding
to the integration constants $u_i (0)$'s,  $0 \le i < k$.  Any two horizontal lifts  of the same
curve differ by a Carnot translation.  See below.  

Our  frame $\{ X_1, X_2 \}$ generates a $k+2$-dimensional nilpotent Lie  algebra $\mathfrak{g}_k$ for which the following commuting relations hold: 
$$X_3 = [X_2, X_1],  X_{4} =[X_{3},  X_1], \ldots , X_{k+2} =[X_{k+1},  X_1],  [X_{k+2},  X_1] = 0,  $$
with
$$X_3 = \dd{}{u_{k-1}},\;\; X_4 = \dd{}{u_{k-2}},\; \ldots\; ,\;\;  X_{k+1} = \dd{}{u_1}, \;\;  X_{k+2} = \dd{}{u_0}.$$
All other Lie brackets $[X_i, X_j]$, $i, j > 1$ are zero.   This algebra is graded nilpotent: 
$$\mathfrak{g}_k = V_1 \oplus V_2 \oplus \ldots V_{k+1},  \;\;\; V_1  = span\{ X_1, X_2\} , \;\;\;  V_i =span\{ X_{i+1}\} ,\;\;  2 \le  i \leq k+1,$$
meaning that  $[V_i, V_j] \subset V_{i+j}$.  (Indeed $[V_1, V_j] = V_{1+j}$,  and $[V_i, V_j] = 0$ if $i, j > 1$).
Thus  $\mathfrak{g}_k$ forms  a $(k+2)$-dimensional graded nilpotent Lie algebra.  The simply connected
Lie group  $G$ associated to any such algebra $\mathfrak{g}$ is, by definition, a Carnot group.  The
exponential map $\mathfrak{g} \to G$ is a diffeomorphism and 
provides $G$ with global coordinates under which the original vector fields are left-invariant
and the multiplication is a `graded polynomial' perturbation of vector addition, with the origin
as the identity.   
Putting a Euclidean structure on its generating level 1 block $V_1$ induces a \sR structure on
the group, with distribution $D$ identified with $V_1$ left-translated about the group. 
  It is   this left-invariant \sR structure which is typically studied
when discussing   Carnot groups.

\subsection{Geodesic equations}  Here is the advertised procedure ( \cite{Monroy1,Monroy2,Monroy3})
for associating geodesics to polynomials in $x$. 
Let  $F(x)$ be any fixed polynomial  in $x$ of degree $k$ or less.
 Solve: 
 \begin{equation}
\ddot x    = - F(x) F' (x),
\label{F-curve1}
\end{equation}
for $x(t)$, insisting that  $x(t)$ also satisfy  the energy constraint 
 \begin{equation}
 \frac{1}{2} \dot x ^2 + \frac{1}{2} (F(x)) ^2  = \frac{1}{2}.
 \label{F-curve2}
 \end{equation}
The energy constraint  is the arc length parameterization condition in disguise. 
Equation (\ref{F-curve1}) is   Newton's equation for the potential $V(x) = \frac{1}{2} (F(x)) ^2$.
  The left hand side of equation (\ref{F-curve2}) is the    conserved total  energy for this Newton's equation.
  
 Having found  such an  $x(t)$, next solve: 
\begin{equation}
\dot u_k (t)  = F(x(t)), 
\label{eq:u}
\end{equation}
for $u_k (t)$.  
The result is a plane curve $c(t) = (x(t), u_k (t))$, $c: I \to \R^2$.  We can always take
this interval $I$ to be the whole real line $I = \R$.   Horizontally lift this plane curve $c$ to form its horizontal
lift $\gamma: \R \to J^k$ using the $k$ `triangular' ODEs  
  (\ref{eq:horiz}) of horizontal lifting. Due to the initial conditions going in to 
horizontal lift, this is not one curve, but a $k$-parameter affine family of such curves parameterized by, for example $u_i (0)$. 
\begin{backgroundtheorem}
\label{thm:back} (See \cite{Monroy1,Monroy2,Monroy3}.) 
The above prescription yields a geodesic in $J^k$ parameterized by
arclength.  Conversely, any arc-length parameterized geodesic in $J^k$ can be achieved
by this prescription applied to some polynomial $F(x)$ of degree   $k$ or less. 
\end{backgroundtheorem} 
\noindent We give an   alternate proof of this  theorem    in  
Appendix A and a second  alternative proof  makes up the final paragraph of  section \ref{sec: mag space}.

The vector fields $X_2, \ldots , X_{k+2}$ span  a codimension one  Abelian subalgebra $\mathfrak{h}$ of $\mathfrak{g}_k$
and the quotient space of $G = J^k$ by the corresponding Abelian group $H$ can be identified with the $x$-axis.   
Translations by elements of $H$ correspond to translations of the coordinates $u_i$.  Since we have left the initial conditions
of our geodesic equations free,  the geodesics determined by a single polynomial $F(x)$ are determined {\it up to left translation}
by elements of $H$.  To translate a geodesic  in the $x$-direction by an amount $x_0$ 
we must  
translate its polynomial : $F(x) \mapsto F(x - x_0)$.  Translation by $x_0$ corresponds to left multiplication by $exp(x_0 X_1)$.
 
 The vector fields $X_3, X_4, \ldots , X_{k+2}$ form, in turn, an Abelian subalgebra $\mathfrak{k} \subset \mathfrak{h} \subset \g$, one which is, moreover,
 normal, being the commutator algebra $\mathfrak{k} = [\g, \g]$. The quotient of $J^k$ by this group is the Euclidean plane $\R^2$
 and the projection $\pi: J^k \to \R^2$ described above corresponds to the quotient projection. 
 On any Carnot group $G$  the analogous projection $G \to G/[G, G] \cong V_1$ is a a submetry.
 In this context the following    principle is  basic to all the work that follows.  
 \begin{proposition}
 \label{prop: lifts of minimizers}
 If $M \to N$ is a submetry, meaning a submersion between \sR submanifolds
 whose distributions have the same dimension, and with the property that the differential of the projection
 is an isometry between distribution planes, then the  
 horizontal lift of a {\it minimizing geodesic} on $N$  is a minimizing geodesic on $M$.
 \end{proposition} 
 \begin{proof}  Points   on $N$ correspond to the  fibers of $\pi$ upstairs on $M$,
 a horizontal curve minimizes between two points of $N$ if and only if its lift minimizes
 between the fibers upstairs.  In particular, the lift  minimizes the
 \sR distance between any two of its points. 
 \end{proof}
 
 \vskip .2cm 
 We apply this principle, that is, proposition \ref{prop: lifts of minimizers} to our case of $\pi: J^k \to \R^2$
The geodesics on $\R^2$ are all known: they are the lines of the first paragraph of this paper 
and they are all global minimizers.  Thus we  have
 a corresponding family of globally minimizing geodesics on $J^k$, the horizontal lifts of lines in the plane.
 These `horizontal lines''  are precisely the curves of the form $t \mapsto h_0 exp(t Y)$
 where $Y = a X_1 + b X_2$ - the left translates by $h_0 \in J^k$  of one-parameter subgroups lying in  the first level
 $V_1$.  One checks without difficulty that these lines  are in bijection with the constant polynomials
 $F(x) = b$, in which case we take $a^2 + b^2 = 1$.    
 
 We are interested in the non-line geodesics, so those geodesics corresponding to non-constant polynomials.
 Multiply the energy equation (\ref{F-curve2})  by $2$ to get    
\beq
(\frac{dx}{dt})^2 + F(x)^2 = 1,
\label{dominant} 
\eeq
 for  $x$ as a function of arclength $t$. 
 Since $\dot x ^2 \ge 0$ everywhere it follows that   $x(t)$   must travel  within one of the intervals on which $F(x)^2 \le 1$.
   We call these the  ``Hill intervals'' of $F(x)$.  There are at most $k$ such intervals
   since their endpoints must be    solutions of the equation  $F(x)^2 = 1$.
    Once we choose one of these Hill intervals $I \subset \R$  for $x$ to travel in,
    the solution $x(t)$ is unique up to a time translation $x(t) \mapsto x(t-t_0)$.  
    To summarize,  every non-line geodesic is determined,
 up to a Carnot translation fixing the   $x$-axis and a time translation,    by a choice of  a degree $k$ polynomial  $F(x)$ together with  one of its  Hill intervals $I$. The endpoints of $I$ satisfy $F(x) = \pm 1$ and the interior points $x$ satisfy $F(x)^2 < 1$. 
\begin{remark}
Given this bijection between geodesics on $J^k$ and the pairs $(F(x), I)$, in the future we will specify a pair $(F(x), I)$ to define a geodesic $\gamma$.
\end{remark}

 According to basic theory  of one-degree of freedom classical mechanical systems,
 there are  three possibilities for the $x$-curve depending on whether or not the endpoints of its Hill interval are critical points of $F$. 
\begin{itemize} 
\item $x$ is periodic of some period $L$: : $x(t+L) = x(t)$.  In this case neither endpoint of $x$'s Hill interval
is a critical point of $F$. These endpoints are referred to as the turning points of the solution. 

\item $x$ is heteroclinic: $t \mapsto x(t)$ traverses its Hill interval exactly once as $t$ varies over $\R$ and does so in a strictly 
monotone fashion.  As $t \to + \infty$, 
 $x(t)$ limits to one   endpoint of its Hill interval,  while as   $t \to - \infty$ it   limits to the other endpoint. In this case both  endpoints  
 of the Hill interval  are critical points of $F$. The solution has no turning point. 
 
\item $x$ is homoclinic:  $t \mapsto x(t)$ traverses its Hill interval twice while $t$ varies over $\R$.  Thus $x(t)$  limits to the same  endpoint $x_0$ of the interval as $t \to \pm \infty$.
It hits the other endpoint  $x_1$ once, at which instant $\dot x = 0$.  We have $F'(x_0) = 0$ while  $F'(x_1) \ne 0$ where $\{x_0, x_1 \}$ are the endpoints
of the Hill interval $I$.  The solution has a single turning point, $x = x_1$. 
\end{itemize} 

In the heteroclinic case, we add one more dichotomy into the mix. 
\begin{definition} 
A  heteroclinic $x$-curve  with Hill interval $[x_0,x_1]$ is said to be of turn-back type if $F(x_0) \ne F(x_1)$, or equivalently, if  
$F(x_0) F(x_1) = -1$.  Otherwise, we say that  the heteroclinic $x$-curve is of direct type, in which case  $F(x_0) = F(x_1)$, 
or equivalently, if $F(x_0) F(x_1) = +1$. 
\end{definition}

 \begin{definition} A non-line  geodesic  is  called x-periodic,  heteroclinic, or homoclinic  
according to whether its   $x$-curve is periodic, heteroclinic or homoclinic. Similarly, we can speak
of non-line geodesics as being heteroclinic of direct type or of turn-back type. 
\end{definition}

 \subsection{Main Results}

 \begin{theorem}\label{theorem-1}  The following classes of geodesics in $J^k$  fail to be   globally minimizing
 \begin{itemize} 
 \item{(i)} those which are  $x$- periodic.
 \item{(ii)} those which are heteroclinic of  turn-back type.  
 \end{itemize} 
 \end{theorem}
 \noindent This theorem is proved in section \ref{sec: proof of B}. It is perhaps not a big surprise to a few experts who
 can prove both (i) and (ii)   by the means by which we will prove item (ii). 
 
 What remains as possible globally minimizing geodesic candidates   are  homoclinic  geodesics and the heteroclinic geodesics  of direct type.  Ardentov and Sachkov \cite{Ardentov1,Ardentov2}   established the minimality
 of the homoclinic geodesics corresponding to $F(x) = a x^2 - 1,  a > 0$  when $k =2$.   Their work provided much of our inspiration. 
  The plane curve for these geodesics will be called the   Euler kink. (Other names for this plane curve  are   syntactrix and  convict's curve.)  See 
  the middle panel of figure  \ref{fig:kink}. 
  In subsection \ref{more submetries} near the end of this section we observe  that the Euler kink   continues to be globally  minimal for all   $k > 2$.


Theorem \ref{theorem-1} also  implies  that any global minimizer which is not a line  must be bi-asymptotic to `vertical lines',
meaning the horizontal lifts of lines of the form $x = const.$.
This fact    is known
to a few experts who understand the results of  \cite{LeDonne}.  
 
 We proceed to our main new result.  

\begin{definition} 
\label{def: seagull}  Call a real polynomial $F(x)$ a ``seagull polynomial'' if 
it is   even, has maximum   $1$, has $0 < F(0) < 1$
and its only critical points are $0, \pm a$ where $F(\pm a) =1$.   
\end{definition} 
The graph of a  seagull polynomial $F(x)$ is qualitatively that of  a   double well potential, reflected about the $x$-axis. 
Since   $a$ and $-a$ are double roots of  $F(x) = 1$ 
we have   
\beq
\label{eq: W}
1- F(x) = (x^2 - a^2)^2 W(x)
\eeq
 with $W(x) > 0$ and the only critical point of $W(x) = 0$ is $x = 0$.  The
set of   seagull polynomials of degree $2k$ forms  a non-empy  open set
of dimension $k-1$ within the $k+1$ dimensional space of   even polynomials of degree $2k$.
To get this count write   $W(x) = a_0 + a_1 x^2 + a_2 x^4 + \ldots + a_{k-2} x^{2(k-2)}$, 
  insist   that $a_0$ and $a_{k-2}$ are positive, $0 \leq a_i$ for $0 < i <k_2$ and  impose  $a^4 a_0 < 1$, these last conditions imply that the set is open. Take the maximum point at $x = a$ as
as an additional parameter.

\begin{figure}%
       {{\includegraphics[width=2.5cm]{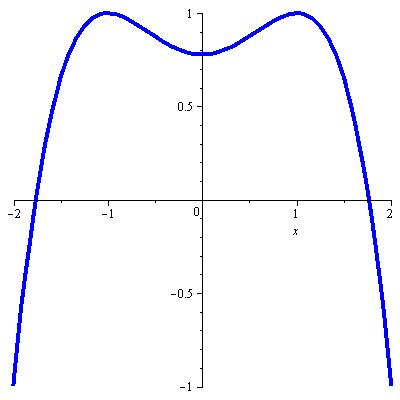} }}%
      {{\includegraphics[width=2.5cm]{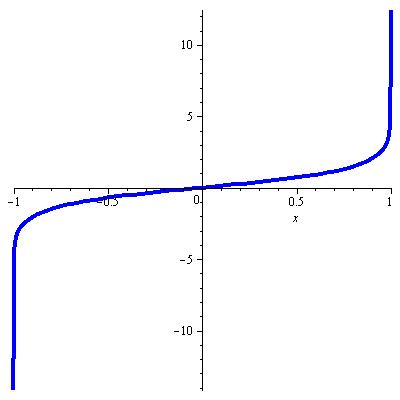} }}
      \caption{ The graph of a seagull polynomial (left panel) and the projection of its
      associated geodesic to the $(x, u_k)$ plane.} 
      \label{fig:seagull} 
\end{figure}

\begin{theorem} \label{main}     
There is a non-empty $8$ dimensional open set of seagull polynomials of degree $18$
 all of which  yield globally minimizing geodesics of heteroclinic type on $J^k$, for any $k \ge 18$.  
 See definition \ref{def: specific class} for specifics regarding this set of polynomials. 
 The Hill intervals $I$ of these geodesics are $[-x_0,x_0]$ where $\pm x_0$ are the global maximum points of the seagull polynomial. (See figure \ref{fig:seagull}.) 
\end{theorem}

The restriction to degree $18$ specified within definition \ref{def: specific class} below occurs at only one step in our proof,
the ``$Leg 3$'' step near the end.  We are confident  the theorem holds
for all even degrees.  
   

\subsection{Miscellany. } 
 
 \subsubsection{ Vertical lines as abnormal geodesics.}  
  The  geodesics for the constant polynomials $F = \pm 1$ form a special class of lines
  called ``vertical lines'' since they are given by   $\dot x = 0$ in   the $(x, u_k)$ plane.  These are   precisely the   abnormal, or {\it singular}  geodesics of $J^k$.   See 
 \cite{tour}, \cite{monster} or \cite{BryantHsu}.   What makes them special in the metric category
 is that they are geodesics, independent of the variable inner product placed on the distribution planes. 
 Theorem \ref{theorem-1} implies  that any non-line  global minimizer is asymptotic to some  vertical line
 as arclength $t \to +\infty$ and to a different vertical line as $t \to -\infty$.  These   distinct lines   have the same projection to the plane 
 in the homoclinic case.   This  fact   instantiates a general  theorem found in \cite{LeDonne}.

 \vskip .3cm

\begin{figure}
     \centering
     \begin{subfigure}[b]{0.27\textwidth}
         \centering
         \includegraphics[width=\textwidth]{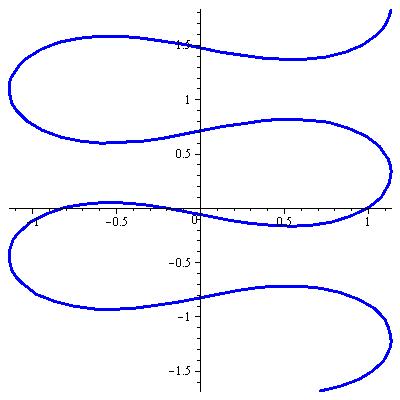}
         \caption{$F(x) = \frac{1}{\sqrt{2}}(2x^2-1)$}
     \end{subfigure}
     \hfill
     \begin{subfigure}[b]{0.27\textwidth}
         \centering
         \includegraphics[width=\textwidth]{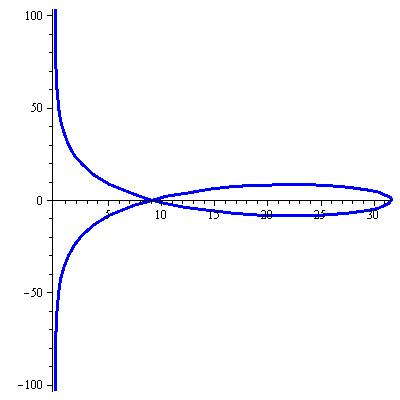}
         \caption{$F(x)=1-\frac{2}{1000}x^2$}
     \end{subfigure}
     \hfill
     \begin{subfigure}[b]{0.27\textwidth}
         \centering
         \includegraphics[width=\textwidth]{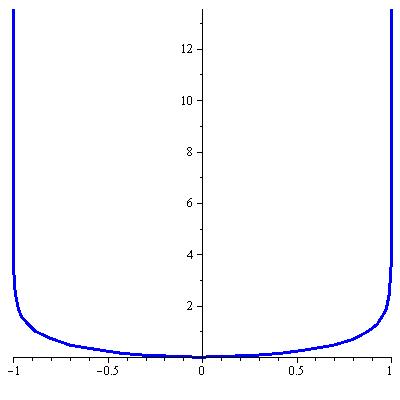}
         \caption{$F(x) = \frac{1}{2}(x^3-3x)$}
     \end{subfigure}
        \caption{The left  panel shows a periodic geodesic.  The middle panel  shows 
    the  kink curve. The right panel shows 
    a planar curve which is the projection of a critical geodesic of turnback type.}
        \label{fig:kink}
\end{figure}

\noindent See figure \ref{fig:kink} for some representative examples.  
%

 \noindent See   the rightmost panel of figure \ref{fig:kink}.  In the turnback type we have $F(x_1) = -F(x_0) = \pm 1$ so that the shape of the plane curve
looks like a giant $U$ according to equation (\ref{eq:u}) with $u_k$ reversing course as $t \to + \infty$ and traveling back the way it came
from in the distant past.  In the direct case, the asymptotic 
direction of motion of $u_k$  is the same in  the distant past and distant future.    
 
\subsubsection{Scaling.} 
\label{sec: scaling} Carnot groups  admit dilations 
$\delta_h: G \to G$, $h \in \R \setminus \{0\}$. The $\delta_h$ comprise a  one-parameter group of  automorphisms of $G$
 which are also  metric dilations: $d(\delta_h g, \delta_h y) = |h| d(g, y)$. 
 If $\gamma(t)$ is a geodesic parameterized by arc-length  then so is $$\gamma_h (t) = \delta_{\frac{1}{h}} \gamma (ht),$$
 for any $h \ne 0$.
 The   Carnot dilation  on $J^k$ is 
$$\delta_{h} (x, u_k, u_{k-1}, \ldots ,u_0) := (h x , h u_k, h ^2 u_{k-1}, h ^3 u_{k-2} ,\ldots, h ^{k+1} u_0).$$
 One verifies by direct computation,  using the geodesic equations
(\ref{F-curve1}), (\ref{eq:horiz}),  and (\ref{eq:u}),  that if $F(x)$ is the polynomial yielding the non-line geodesic $\gamma (t)$ then $F_{h} (x):= F(h x )$
is the polynomial yielding  the scaled non-line geodesic $\gamma_h (t)$.   

 \subsubsection{ The other jet submetries.}
 
 \label{more submetries}
 
 We discussed the submetry $J^k \to \R^2$.  As a metric space, this $\R^2$ coincides with $J^0$. 
 These fit into a family of     \sR submersions to lower level jets,  $\pi_{k,n}: J^k \to J^n$ with   $n < k$, so that
 $\pi = \pi_{k,0}$, 
  $$\pi_{k,1}: (x,u_k,\dots,u_1,u_0) \mapsto  (x,u_k,u_{k-1}).$$
$$ \vdots $$
$$\pi_{k, k-1}: (x,u_k,\dots,u_1,u_0) \mapsto  (x,u_k,\dots,u_2,u_1) : = (x, v_{k-1}, \dots , v_1, v_0).$$
The last map $\pi_{k,k-1}$
realizes the quotient map $J^{k-1} \cong J^k/\exp(V_{k+2})$. 
We   identify this quotient space with $J^{k-1}$ 
by  shifting the  meaning of coordinates - the old $u_0$
 has been  projected out, and its derivative $u_1 = du_0/dx$ is set to $v_0$ which now plays the role of the function whose jet we are taking
 when forming $J^{k-1}$. The old $u_k$ continues its role as the `fiber coordinate' of jet space, but this
 time now  in the role of  the $(k-1)$th derivative of $v_0$ with respect to $x$.  
 
 It  follows from the basic principle, proposition \ref{prop: lifts of minimizers},  that a globally minimizing geodesic corresponding to some
 degree $n$ polynomial $F(x)$ persists by horizontal lift   to yield a globally minimizing geodesic for all higher $k > n$.  
 In particular the  geodesic which projects to the Euler kink  continues to be a global minimizer for all $k > 2$.

\section{Proof of  theorem \ref{theorem-1} }   
\label{sec: proof of B}

\subsection{Case (i) the x-periodic case} 

\begin{proposition} Let $K$ be the following vector field 
$$K = \sum_{i=0}^k \frac{x^{k-i}}{(k-i)!} \frac{d}{du_{i}}, $$
then $K$ is a Killing vector field. 
\end{proposition}
\begin{proof}
First let us introduce a equivalence definition for a Killing vector field. Let $P_1,P_2 : T^*J^k \to \R$ be the momentum functions of the vector fields $X_1,X_2$, see \cite{tour} 8 pg. In terms of traditional cotangent coordinates $(x,u_k,\cdots,u_0,p_x,p_{u_k},\cdots,p_{u_0})$ for $T^*J^k$, we have
$$ P_1 = p_x + \sum_{j=0}^{k-1} u_{j+1} p_{u_{j}}, \;\; P_2 = p_{u_k}. $$
Then the Hamiltonian governing the geodesic on $J^k$ flow is $H = 1/2(P_1^2+P_2^2)$. So $K$ is a Killing vector field if and only if its momentum function $P_K$ Poisson commute with $H$, that is $\{P_K,H\} = 0$. Then it is enough to prove that $\{P_K,P_1\} = 0$ and $\{P_K,P_2\} = 0$, when
$$ P_K = \sum_{i=0}^k \frac{x^{k-i}}{(k-i)!} p_{u_{i}}. $$
We have that $\{P_K,P_2\} = 0$, since $P_K$ does not depend on $u_k$. Then we will focus on the first bracket, 
\begin{equation*}
\begin{split}
\{P_K,P_1\} &= \{p_x,  \sum_{i=0}^k \frac{x^{k-i}}{(k-i)!} p_{u_{i}} \} + \{\sum_{j=0}^{k-1} u_{j+1} p_{u_{j}} , \sum_{i=0}^k \frac{x^{k-i}}{(k-i)!} p_{u_{i}}\} \\
            &= \sum_{i=0}^k p_{u_{i}} \{p_x,  \frac{x^{k-i}}{(k-i)!} \} + \sum_{j=0}^{k-1} \sum_{i=0}^k p_{u_{j}}  \frac{x^{k-i}}{(k-i)!} \{ u_{j+1} ,  p_{u_{i}}\} \\
            &= - \sum_{i=0}^{k-1} p_{u_{i}} \frac{x^{k-i-1}}{(k-i-1)!} + \sum_{i=0}^{k-1} \sum_{j=0}^k p_{u_{i}}  \frac{x^{k-j}}{(k-j)!} \delta_{i+1,j},   \\
\end{split}
\end{equation*}
where we used in the last line that the term from the first sum when $i = k$ does not depend on $x$, then we can sum until $i =k-1$. In the second sum we switch the place of $i$ and $j$, also we used $\delta_{j+1,k-i}$ that is the Kronecker delta, $\delta_{i+1,k-j}$ is equal to $1$ when $i+1 = j$ and zero otherwise, then 
\begin{equation*}
\begin{split}
\{P_K,P_1\} &=  - \sum_{i=0}^{k-1} p_{u_{i}} \frac{x^{k-i-1}}{(k-i-1)!} + \sum_{i=0}^{k-1} p_{u_{i}}  \frac{x^{k-i-1}}{(k-i-1)!} = 0.   \\
\end{split}
\end{equation*}
\end{proof}

The last proposition implies that the flow of $K$ generates a \sR isometry. Now we are ready to prove case (i) from \ref{theorem-1}.

\begin{proof} Case (i) from \ref{theorem-1}:
Let $\gamma$ be a geodesic  for the polynomial  $F(x)$ whose  $x$-curve  $x(s)$ is periodic of period $L$.  
Let  $[x_0,x_1]$ be the   Hill interval for $x(s)$  so that  $x_0, x_1$ are turning  points for $x$.  
By performing an $x$-translation we may assume that $x_0 = 0$ and by an $s$-translation
that $x(0) = 0 = x_0$.  Then $x(L/2) = x_1$ and $x(L) =0$.   We claim that $\gamma(L)$ is conjugate
to $\gamma(0)$ along $\gamma$.

Next,  observe that at the  turning points $s =0, L/2, L , 3L/2, 2L, \ldots$  of the $x$-curve we have  $\dot x (s) = 0$ so that $\gamma$ is tangent to the vertical direction $\dd{}{u_k}$ at these times.  In particular, by reversing directions if necessary, we have that $\dot \gamma(0) = \dot \gamma (L)  =  \dd{}{u_k}$.
Now consider the following
two Jacobi fields for $\gamma$:
\begin{equation*}
\begin{split}
W_1(s) &= K \text{ restricted to } \gamma, \\
W_2 (s) &= \dot \gamma (s).
\end{split}
\end{equation*} 
 Since $x(kL) = 0$ we have that $W_1 (jL) = \dd{}{u_k}$, $j =0, 1, 2, \ldots$
 so that   $W_1(0) = W_2 (0) = W_1 (L) = W_2 (L)$.
Since  the space of Jacobi fields is a linear space
 so that $J : = W_1 -W_2$ is again a  Jacobi field for $\gamma$ and this field now
  vanishes at every $s = jL, j= 0, 1, \ldots$.
 In the interior of the interval $(0,L)$ the field $J$ is  not identically zero since  $\dot x(s) \ne 0$ for $0 < s < L/2$. 
 It follows that $J$  contributes at least $1$ to the nullity of the Hessian of the action - so   the squared length 
 functional $\int_0 ^s \frac{1}{2} \| \dot \gamma (s)\|^2 ds$ - thus  establishing that the times   $s = L, 2L, \ldots $ are conjugate times  to $s =0$ 
 along $\gamma$. It   follows by  standard calculus of variations that  the geodesic
 $\gamma$  fails to minimize   beyond $L$.  
\end{proof}

\vskip .3cm

\subsection{ Case (ii) : the heteroclinic turn-back case.} 
 
Our proof relies on  the method of blowing-down geodesics as  explained by   Hakavuouri-Le Donne \cite{LeDonne}. Suppose that $\gamma: \R \to G$ is a rectifiable curve  in a Carnot group $G$.
For $h \in \R^+$ form
$$\gamma_h (t) = \delta_{\frac{1}{h}} \gamma (ht),$$
where $\delta_h: G \to G$ is the Carnot dilation.    One easily checks that
if $\gamma$ is a geodesic then so is $\gamma_h$ for any $h > 0$.

\begin{definition}  A {\it blow-down} of $\gamma$ is any limit curve
$\tilde \gamma = \lim_{k \to \infty} \gamma_{h_k}$
where $h_k \in \R$ is any sequence of scales tending to infinity with $k$,
and the limit being uniform on compact sub-intervals. 
\end{definition}

Hakavouri  and LeDonne \cite{LeDonne}   prove the following lemma:
\begin{lemma}  If $\gamma$ is globally minimizing geodesic parameterized by arclength
then any one  of its blow-downs $\tilde \gamma$ is also a globally minimizing geodesic parameterized by arc-length.
\end{lemma} 

\begin{proof} Proof of case (ii) from \ref{theorem-1}:
   The projection $\pi$  to the $(x, u_k)$ plane  of a heteroclinic turnback  geodesic $\gamma$ 
lies between two vertical  lines $x = x_0$ and $x = x_1$ and its height $u_k$ achieves   a  global maximum or  minimum
at some point $P$  in between these   lines. (See the right panel of figure \ref{fig:kink}.)  The geodesic  $\gamma$ is asymptotic to  
of one of these  vertical lines as $t \to -\infty$ and to the other as $t \to + \infty$. 
 Using a time translation if needed, we may  assume that the extremal point  $P$   occurs when 
  $t = 0$ and by using a translation we can assure that  $P =0 \in J^k$. 
And by using the dilation $\delta_{-1}$   we can  
assure that $P$ is a global minimum point  for $u_k$, so that   $u_k  >  0$ everywhere else along the curve.
Let $[x_0,x_1]$ be the geodesic's Hill interval. Then, upon dilating we have that
  $\pi \circ \gamma_h$ lies between the  two vertical lines
$x = x_0/h$ and $y = x_1/h$ and is asymptotic to one of them as $t \to -\infty$ and the other as $t \to + \infty$,
while for all $\gamma_h (0) =0$.  It follows that any blow-down of $\gamma$ consists of the horizontal lift of the vertical ray
$x = 0,  u_k \geq 0$, traversed twice, once coming in from infinity, hitting zero, then reversing course and heading back
out to infinity.  But a ray, traversed twice is not a minimizing geodesic, since it is not smooth curve on the $x-u_k$ plane (any geodesic on a \Ri manifold is smooth), and neither are any of its horizontal lifts.
So $\gamma$ cannot itself cannot be a globally minimizing  geodesic, for if it were, the Hakavuori-Le Donne  lemma would imply that the vertical
ray, traversed there and back, is a globally minimizing geodesic.  But a curve which retraces its own path is never a minimizing geodesic: simply chop off and shorten the path by stopping before the turn-around point $P$ and turn back earlier. 
\end{proof}

\section{Setting up for theorem \ref{main}.}\label{sec: magnetic}

\subsection{The intermediate Magnetic space}

   As a  first  step towards proving theorem \ref{main}    we  factor the  \sR submersion   $\pi: J^k \to \R^2$ into the  product of two
\sR submersions:
\begin{equation}
\label{factor}
\pi = pr \circ \pi_F,
\end{equation}
where the target of $\pi_F$ is an   intermediate 3-dimensional \sR space denoted by 
 $\R^3 _F$ whose geometry  depends on the choice of polynomial $F(x)$ and which we refer to 
 as a `magnetic \sR structure'.   Thus we will have  \sR submersions
 $$ \xymatrix{J^k \ar[r]^{\pi_F} & \R^3 _F  \ar[r]^{pr} & \R^2 _{x , u_k}  }, $$
if $x, y, z$ are coordinates on this intermediate space
then
the distribution $D_F$ on $\R^3 _F$ is defined  by the single Pfaffian equation:  
 \begin{equation} 
 D_F:  dz - F(x) dy = 0, 
 \label{interm_distrib1}
  \end{equation} 
while the metric $ds^2$  on the two-planes $D_F$ is defined  by  
  \begin{equation} 
ds^2 =  (dx^2 + dy^2)  |_{D_F}, 
 \label{interm_distrib2}
\end{equation} 
 and  the projection to the plane $\R^2$ by
$$pr(x, y, z) = (z, y) = (x, u_k).$$

Before we describe the projection we pause to explain why we've used the term ``magnetic''.
The motion of a  particle of charge $e$   moving non-relativistically 
 in the Euclidean plane under the influence of a magnetic field of strength $B(x,y)$
orthogonal to the plane is given by $\ddot c = e B(c) {\mathbb J} \dot c$. 
Let $A = A_1 (x, y) dx + A_2 (x,y)dy$ be a vector potential for $B$, meaning that $dA = B dx \wedge dy$.
The Hamiltonian system on $T^* \R^2$ having   Hamiltonian $H = \frac{1}{2} (p_x - e A_1 (x, y))^2 + (p_y - e A_2 (x, y))^2$
generates the motion of this particle. (See, eg \cite{Landau}.) 
Introduce a third variable $z$ with conjugate momentum $p_z$ so that $H$ becomes
  $H = \frac{1}{2} (p_x - p_z A_1 (x, y))^2 + (p_y - p_z A_2 (x, y))^2$ 
on $T^* \R ^3$.  This is the \sR kinetic energy for the \sR structure  on $\R^3$
defined by the distribution $D = ker(dz - A)$ 
with inner product $dx^2 + dy^2 |_D$.  Since $H$ is independent of $z$ we 
have that $p_z$ is constant along trajectories and we identify this constant with the charge $e$.    (See, eg \cite{tour}.) We call any \sR structure of this form
on $\R^3$ a magnetic \sR structure.    Our $\R^3 _F$  is such a structure with $A = F(x) dy$.

To construct the projection  $\pi_F$  
expand out 
$$F(x) = \Sigma \frac{a_j}{j!} x^j $$
and use the {\bf alternate} coordinates $\theta_j$ for   $J^k$ as those described in \cite{Monroy1} and \cite{ Monroy2} with an index swapping and sign change:    
$$\theta_0 = u_k, \cdots ,\; \theta_j = \sum_{i=0}^{j} (-1)^i\frac{x^{j-i}}{(j-i)!}u_{k-i} ,  \cdots,  \;\theta_k = \sum_{i=0}^{k} (-1)^i\frac{x^{k-i}}{(k-i)!}u_{k-i} , $$
(The $u_j$ of  Monroy-Perez and Anzaldo Meneses are not exactly the same as ours,
rather they are related to ours by an index swapping and a sign change. To be precise,  that is, if  $v_j$ denotes the $u_j$ in \cite{Monroy1}  
then $v_j = (-1)^j u_{k-j}$.)

Written in $\theta$-coordinates our frame for the distribution $D$ on $J^k$  is  
$$X_1 = \dd{}{x},\qquad X_2  = \dd{}{y} + \sum  \frac{x^j}{j!} \dd{}{\theta_j}, $$
where $y = u_k = \theta_0.$

Define 
$\pi_F: J^k \to \R^3$
 to be    linear in the $\theta$ coordinates and with  coefficients constructed from the  
 scaled coefficients   $a_j$ of  our  chosen $F(x)$: 
\begin{equation*}
\begin{aligned} 
\pi_F (x, \theta_0 , \theta_1, \ldots , \theta_{k}) = (x,\theta_0, \sum_{i=0} ^k  a_i \theta_i): = (x,y, z) \\
\end{aligned}
\end{equation*}
in the original coordinates the projection looks as
\begin{equation*}
\begin{split}
\pi_F (x, u_k , u_{k-1}, \ldots , u_{0}) = (x,u_k, \sum_{i=0}^{k} (-1)^{k-i} u_{i} \frac{d^{k-i}F}{dx^{k-i}}(x))  .  \\
\end{split}
\end{equation*}
From the  linearity of $\pi_F$ in these coordinates  we easily compute   
$$\pi_{F *} X_1 = \dd{}{x} \;\; \text{and}\;\; \pi_{F *} X_2 = \dd{}{y} + F(x) \dd{}{z},$$
which is an \on  horizontal frame for the \sR structure $\R^3 _F$ as defined by    equations (\ref{interm_distrib1}) and  (\ref{interm_distrib2}), 
establishing the  factorization (\ref{factor}).

\subsection{ Magnetic geodesics.} 
\label{sec: mag space}
The \sR kinetic energy Hamiltonian for $\R^3 _F$ is 
$$H_F = \frac{1}{2}p_x ^2  + \frac{1}{2} (p_y + F(x) p_z)^2 .$$
The projection   $(x(s), y(x), z(s))$  to    $\R^3 _F$ of any solution  curve  
$$(x(s), y(x), z(s), p_x (s), p_y(s), p_z (s)) \in T^* \R^3 _F ,$$  to Hamilton's equations
for $H_F$  is a geodesic.  The geodesic is parameterized by arc-length if   $H_F = 1/2$.
Since no $y$'s or $z$'s occur in $H_F$, the time derivatives of the  momenta $p_y , p_z$ are zero so
we have that 
$$p_y = a , \;\;\text{and}\;\; p_z = b $$
with $a$ and $b$ constant.

It makes sense to call  $p_x = p $  momentum for $x$, 
we then have  
$$H_F = \frac{1}{2} p ^2 + V(x), \qquad V(x) = \frac{1}{2}(a + b F(x)) ^2. $$
Hamilton's equations for the pair $(x, p)$ are identical to the 1st geodesic equation,
equations (\ref{F-curve1}) and (\ref{F-curve2}) upon replacing  $F$  by $G = a + bF$. Thus the $x$-component of our
geodesic is, following definition \ref{F-curve1}, an $x$-curve for $G$.  
\begin{definition} The {\it pencil} of the polynomial $F$ is the two-dimensional linear space
of polynomials having  the form
\begin{equation}
G(x):= a + b F(x),
\label{eq: pencil}
\end{equation}
where $a, b$ are arbitrary real constants. 
\end{definition}
 What about the   remaining components $y(s), z(s)$ of our geodesic?  
 Using that $H_F = \frac{1}{2} (p^2 + P_y ^2)$ with  
 $P_y = p_y + p_z F(x) = a + b F(x) = G(x)$
 and writing out Hamilton's equations  for $y$ and $z$ we get
\begin{equation}\label{F-magnetic}
\begin{split}
\dot y & = G(x),        \\
\dot z & = G(x)  F(x). \\
\end{split}
\end{equation}
The first equation is  equation (\ref{eq:u}) after replacing  $F$  by $G = a + bF$.
The second equation, upon substituting in the first, says that
$\dot z - F(x) \dot y = 0$, which simply says that $(x(s), y(s), z(s))$ is the horizontal
lift of the plane curve 
$$(x(s), y(s)) = (x(s), u_k (s)).$$

Since $\pi_F$ is a \sR submersion, the horizontal lift of any  $\R^3_F$-geodesic for $F$ is a
geodesic on $J^k$.   Horizontal lift is given by lifting the plane curve $(x(s), y(s)) = (x(s), u_k (s))$,
which is to say, by  equations (\ref{eq:horiz}). We have proven
\begin{lemma} Every $\R^3 _F$ geodesic is the $\pi_F$-projection of a geodesic  in $J^k$
corresponding to some $G$ in the pencil of $F$.  Conversely, the  horizontal lifts to $J^k$ 
of $\R^3_F$-geodesics are precisely those  geodesics  
corresponding to   polynomials  in $F$'s pencil.  
\end{lemma} 

{\bf Remark. }
As an immediate corollary to the lemma we get:

{\sc Proof of  Theorem \ref{thm:back} (``Background Theorem'').} 
Take $a =0, b =1$ so that  $G = F$.  The lift of the geodesic in $\R^3 _F$ corrresponding to $G= G$
is a $J^k$-geodesic.

\subsection{Periods I}
\label{sec: Periods} 
\vskip .2cm

\begin{proposition}\label{y-z as x function} Let $c(t) = (x(t), y(t), z(t))$ in $\R^3 _F$ be the projection of an $x$-periodic geodesic in $J^k$
having  corresponding polynomial   $G = a + b F$ and Hill interval $[x_0, x_1]$. (Recall that $G(x_0)^2 = G(x_1)^2 = 1$.) 
Then the  $x$-period   is  
\beq
\label{eq: period} L  = 2  \int_{x_0}^{x_1} \frac{dx}{\sqrt{1-G^2(x)}},
\eeq
and is twice the time it takes the $x$-curve to cross its Hill interval exactly once. 
After one period   the changes $\Delta y = y(t_0 +L) - y(t_0)$ and $\Delta z = z(t_0 +L) - z(t_0)$ undergone by $y$ and $z$ are given by  
\[   \Delta y  = 2 \int_{x_0}^{x_1} \frac{G(x)dx}{\sqrt{1-G^2(x)}} ,\; \Delta z  =2  \int_{x_0}^{x_1} \frac{F(x)G(x)dx}{\sqrt{1-G^2(x)}}.  \]
\end{proposition}

 \begin{proof}  
Along any arc of $c$ for which  $x(t)$ 
is monotonic we can re-express the   curve
as a function of $x$
instead of $t$   using $(dx/dt)^2 = 1 - G(x)^2$ or $\frac{dx}{dt} = \pm \sqrt{ 1 - G(x)^2}$.   
Note the sign of the $\pm$ changes each time  $x(t)$ reflects off of an endpoint of its Hill interval.
Then, as is standard in mechanics, the total period $L$  is twice  the time $\Delta t$ required to cross the Hill interval.  
We have 
$$\Delta t= \int dt = \int \frac{dt}{dx} dx = \int_{x_0} ^{x_1}  \frac{1 }{\sqrt{1-G^2(x)}}dx.$$

For the other two periods use that 
 the  differential equations (\ref{F-magnetic})  assert that $dy = G(x(t))dt$ and $dz=F(x(t))G(x(t)) dt$. Choose $x(t)$ so that $x(0) = x_0$.  Then  $x(L/2) = x_1$.
 and  $x(L/2 + t) = x(L/2 -t)$.  It follows that the change in $y$ and $z$  over a full period is   twice their change
 over a half-period.    Since $dt = \frac{1 }{\sqrt{1-G^2(x)}}dx$ on the  first period  
 we get the result  provided we start at $x_0$ at time  $t = 0$.  To see that the result of $\Delta y$  is independent of
 the starting point, differentiate $y(t + L) - y(t)$ with respect to $t$.  The derivative is $G(x(t+L)) - G(x(t))$. But $x(t+L) = x(t)$ 
 so this derivative is zero.  The same proof works for $\Delta z$.   
\end{proof}   
  
{\bf Remark.} The above equation for $L$ is the well known formula for the period of a one-freedom degree mechanical system. Where $x_0$ and $x_1$ are the points when the potential energy is equal to total energy of the system, $F(x)^2 = 1$, and we use the reversibility of the system, if $x(t)$ is a solution to the Newton's equation then $x(-t)$ is too, to assure that the time that takes to the particle to travel from $x_0$ to $x_1$ is equal to the one from $x_1$ to $x_0$, then the period $L$ is equal to the time that takes to the particle to come back to the initial point and is independent of the initial point and initial time.     

{\bf Remark.} $\Delta y$ and $\Delta z$ are also independent of $t_0$, which is to say, of  the    initial point of the curve.

\section{Calibrations: a Hamilton-Jacobi method for local minimality}\label{sec-H-J}

Suppose $H: T^*Q \to \R$ is a Hamiltonian on some standard phase space $T^*Q$.
The  associated time-dependent Hamilton-Jacobi  equation is the PDE 
$$H(q, dS(q)) = const,$$ 
to be solved for  a function $S:Q \to \R$.   
For lines in  Euclidean geometry, we take $Q = \R^n$,   $H(q,p) = \frac{1}{2}\|p \|^2$,    the Hamiltonian of
a free particle, and    the constant to be $1/2$. Then the Hamilton-Jacobi   PDE  reads $\|\nabla S \| = 1$
and in this guise is often  called the {\it eikonal equation} -- the equation of light rays. 
The integral curves of the gradient flow,  $\dot q = \nabla S (q) $,  are straight lines. A typical solution
$S$ has the form  $S(q) = dist(q, C)$ where $C \subset \R^n$ is a closed set. 
All of this extends to the Hamilton-Jacobi equation associated to the geodesic
equations in  \Ri and in \sR geometry. 
 See pp 14-15 of \cite{tour} for details.

In the \sR case the    Hamilton-Jacobi equation associated to geodesic flow   is called the ``eikonal equation'' following the usage of geometric optics.  This equation reads  
\beq
\| \nabla_{hor} S \| =1,
\label{eq: eikonal} 
\eeq
where $\nabla_{hor}$ is the horizontal derivative of $S: Q \to \R$, that is, $\nabla_{hor}S$ is the unique {\it horizontal} vector field satisfying, for every $q$ in $Q$,
$$ \langle \nabla_{hor}S,v\rangle_q = dS(v), \;\;\text{for every} \;v \in D_q . $$
Where $\langle \; , \; \rangle$ is the \sR inner product.

A more careful definition of ``globally minimize''  is in order
\begin{definition} 
\label{def: global min2} 
A. Let $\Omega \subset Q$ be a domain within a \sR manifold $Q$ and $I \subset \R$  a closed bounded interval.
We say that 
 $c: I \to \Omega$  globally minimizes within $\Omega$   if  whenever $\tilde c :J \to \Omega$ 
 is any  smooth horizontal curve   lying in $\Omega$ and  sharing      endpoints with $c$, then  $\ell(c) \le \ell(\tilde c)$. 

B. If  the interval $I \subset \R$ is  not closed and bounded  then we say that $c: I \to \Omega$ is  ``globally minimizing within $\Omega$'' if every   closed bounded sub-arc
$c([t_0, t_1])$ of $c$ , $t_0, t_1 \in I$  is globally minimizing within $\Omega$ in the sense of A.  
\end{definition} 

\begin{remark}
In part A of the definition we could have replaced ``any smooth horizontal curve lying in $\Omega$ by  ``any continuous curve lying in $\Omega$''
without changing the meaning.  The reason is that  the subRiemannian length functional 
$c \mapsto \ell(c)$ satisfies  the property that if $\ell(c) < \infty$ then  $\ell(c) = \lim \ell(C_i)$ where the $C_i \to c$ is any sequence of  smooth
paths  $C_i$ which converge to $c$ in  either  the $H^1$  or the Lipshitz sense. See \cite{tour} for details. 
\end{remark}

\begin{proposition}
If $S$ is a  $C^2$ solution of the eikonal equation (\ref{eq: eikonal})  defined  on a simply connected domain $\Omega \subset Q$,
then the integral curves of its   horizontal gradient flow 
$\dot c = \nabla_{hor} S (c)$ are \sR  geodesics which {\it globally minimize}  within the domain $\Omega$.  
\label{HJ_min} 
\end{proposition} 
\begin{proof} 
 Let $A, B$ be the shared endpoints of our competing curves $c, \tilde c$. Then Stokes's theorem imply:
$$\int_c dS = \int_{\tilde c} dS = S(B) -S(A).$$
But for any smooth curve $\gamma$ in $\Omega$ we have that 
$$\int_{\gamma} dS = \int \langle \nabla S, \dot \gamma \rangle dt
\le \int_{\gamma} \| \dot \gamma \|  \| \nabla_{hor} S \| dt = \int_{\gamma} \| \dot \gamma \| dt = \ell(\gamma).$$
Equality holds in this series of inequalities if and only if $\dot \gamma = f \nabla_{hor } S$ for some positive scalar $f$,  that is, if and only if $\gamma$
is a reparameterization of an integral curve of $\nabla_{hor} S$.  Our curve $c$ is such an integral curve, that is, $\dot{c}(t) = (\nabla_{hor} S)_{c(t)}$, then 
$$dS(\dot c) = \langle \nabla_{hor} S, \dot{c} \rangle = \langle \nabla_{hor} S, \nabla_{hor} S \rangle =
1 .$$
 Any other competing curve lying in $\Omega$ satisfies 
 $$dS(\dot {\tilde c}) =  \langle \nabla_{hor} S, \dot {\tilde c} \rangle < \| \dot {\tilde c} \|$$ on an open set of points, it is a strictly inequality since $\tilde{c}$ is different of $c$ at least on an open set, so  the above equality becomes
\begin{equation} 
\label{HJineq} 
\ell(c) = S(B)- S(A) < \ell (\tilde c) ,
\end{equation}
where $\ell$ is the \sR length.  
\end{proof}

In the particular case where $Q = \R^3 _F$
we can  simplify the eikonal equation.  Recall the \sR structure on   $\R^3 _F$. (See equations (\ref{interm_distrib1}, \ref{interm_distrib2}).) 
Also recall that we are now denoting our  coordinates on $\R^3 _F$ by  $x,y,z$.  (See the beginning of subsection \ref{subsec: Periods}.) 
Take any $S = S(x, y, z)$,  compute 
$$dS = \dd{S}{x} dx + \dd{S}{y} dy + \dd{S}{z}dz,$$
use that $dz = F(x)dy$ on horizontal planes to see that
$$dS|_D=  \dd{S}{x} dx |_D + ( \dd{S}{y}  + \dd{S}{z} F (x)) dy|_D.$$
Since $dx, dy$ form  an \on coframe for $D^*$ we have  
\begin{equation}
\nabla_{hor} S = \dd{S}{x} E_1 + ( \dd{S}{y}  + \dd{S}{z} F (x)) E_2,
\label{horiz1}
\end{equation} 
where
$$E_1 = \dd{}{x},  E_2 = \dd{}{y} + F(x) \dd{}{z},$$
is the \on frame dual to $dx, dy$.  The eikonal equation for $S$ then reads
\begin{equation}
\label{eik}
( \dd{S}{x})^2 + (( \dd{S}{y}  + \dd{S}{z} F (x))^2 = 1.
\end{equation}

Take the ansatz
\begin{equation}
\label{ansatz}
S(x,y, z) = b z  + a  y + f(x),
\end{equation}
 for a solution $S$ to be associated to the polynomial $G = a + b F$ in the pencil of $F$.  
 Then eq (\ref{eik}) becomes
\begin{equation}
\label{eik2}
f'(x) ^2 + (a   +b F (x))^2 = 1,
\end{equation}
and the associated horizontal gradient flow vector field is
\begin{equation}
\nabla_{hor} S = f'(x)  E_1 + ( a   + b F (x)) E_2.
\label{horiz2}
\end{equation} 

Compare these equations with that of the geodesic equations
for $G$.  The $x$-curve for $G$  satisfies the energy equation:  
$$\dot x ^2 + (a + b F(x))^2 = 1,$$
while the   $u_k = y$ equation for $G$'s geodesic  is  
$$\dot y = a + b F(x).$$
 Solve the energy equation to get
$$\dot x = \pm \sqrt{ 1 - (a + b F (x))^2 }.  $$
Conclude that   the horizontal gradient flow  equation (\ref{horiz2}),  
and our geodesic equations are   identical provided  
\begin{equation}
f '(x) =  \pm \sqrt{ 1 - (a + b F(x))^2 }.
\label{solve}
\end{equation} 
Note  that  the Hamilton-Jacobi equation   (\ref{eik2}) is equivalent to  equation (\ref{solve}),
up to a choice of sign; when  $G = F$ we have  the two solutions
\begin{equation}
S= S (x,z) = \pm \int ^x  \sqrt{1-F(u)^2}du + z,
\label{solve F}
\end{equation}
choose one,  say  the one with $+$ sign.

We now analyze the maximum domains of definition  $\Omega$ for our  $S$. This domain   must exclude points where $1- F(x)^2 < 0$ 
in order for the square root in the integral not to be imaginary.   We also want proposition \ref{HJ_min}
to hold, which requires  that $S$   be $C^2$ on its domain. 
If $[x_0,x_1]$ is one of the Hill intervals and $F'(x_1) \ne 0$ then
$\dd{S}{x} = \sqrt{1 - F(x)^2 }$ fails to be $C^1$ at $x_1$, and moreover $1- F(x)^2 < 0$ for $x = x_1 + \eps$.
In this case we must exclude the plane $x = x_1$ and nearby points with  $x = x_1 + \eps$, $\eps > 0$ from $\Omega$.
Of course a similar argument holds if it is $F'(x_0)$ that is nonzero.  That is: we must exclude   points $(x,y, z)$ for which  $x$ is a  non-critical endpoint of a Hill interval for $F$.  
   It follows that if we want  to use  proposition \ref{HJ_min} on $\Omega$  and if  the $x$-curve   for   $F$   is periodic then we must take $\Omega$ to be the pre-image of the open interval $(x_0,x_1)$ under the x-projection,   excluding   the turning points $x= x_1$ and $x = x_0$ associated to our $x$-curve.
On the other hand, if  $x_1$ is a local maximum of $F(x)^2$ then a Taylor series analysis 
shows that $\dd{S}{x} $ is $C^1$ at $x = x_1$. In this case we can adjoin $x = x_1$ to $\Omega$ and at least a small neighborhood of 
points with $x = x_1 + \eps$, $\eps > 0$.  Indeed we can continue through to the entire neighboring Hill interval $[x_1, x_3)$
adding its pre-image to the domain $\Omega$  of $S$.  In this way we get a  larger domain whose $x$ projection is  $(x_0,x_3)$.  If $x_3$ is again a local maximum for $F(x)^2$
we can continue this process.  Eventually we arrive at  the maximal  domain $\Omega$ for  $S$, a domain of the  form:   
\beq \Omega: = \{ (x, y,z): x \in (\alpha , \beta) \} \subset \R^3 _F,
\label{eq: slab}
\eeq
where $ (\alpha , \beta) = (x_0, x_1] \cup   [x_1, x_2] \cup \ldots [x_{k-1} , x_k)  $ and where each   $[x_i, x_{i+1}]$ is a Hill interval for $F$.
All of the $x_i$ but the endpoints $x_0 = \alpha$ and $x_k = \beta$  
are  local  maxima for $F(x)^2$ having value $F(x_i)^2 =1$.  On the other hand  $x_0 = \alpha$ and $x_k = \beta$ are not local  maxima of $F(x)^2$.

  \begin{definition} We will call  $\Omega$ as described in equation (\ref{eq: slab}) and the description following (\ref{eq: slab})
a slab domain for $F$ associated to any one of the    Hill intervals $ [x_i, x_{i+1}]$, whose interior is  contained by $(\alpha, \beta)$. 
\end{definition}
We have proven: 
\begin{proposition}\label{HJminimality} Let $F$ be a non-constant polynomial. Let  $\gamma:\R  \to J^k$ be  a  
 geodesic for $F$ whose  $x$-curve $x(t)$ has Hill interval $[x_0,x_1]$. Let $c = \pi_F \circ \gamma: \R \to \R^3 _F$. Let $\Omega$ be the slab domain for $F$ 
 associated to   $[x_0,x_1]$ as per the above definition.  Let $I \subset \R$ be an open interval, possibly infinite, possibly all of $\R$,  on which $x(t)$ is strictly monotonic
 and satisfies $x(I) = (x_0,x_1)$, where we have explicitly excluded the case of $x_0 \in x(I)$ and
 $x_1 \in x(I)$.  Then $c|_I$ is a global minimizer within $\Omega$. 
 \end{proposition}

\begin{proof}    
Over the interval $I$ the sign of $\dot x$ is fixed, either plus or minus.  
Choose the sign of the square root in equation (\ref{solve F}) accordingly.  We get a smooth solution $S(x,z) = S(x,y,z)$
to the Hamilton-Jacobi equation and the open arc of our geodesic $c(I)$
is an integral curve of $\nabla_{hor} S$.  Thus $c|_I$ is a global minimizer within
$\Omega$ by proposition \ref{HJ_min}.  
\end{proof}

\vskip .3cm 

It is worth seeing how this argument looks in each of our three cases. 

{\sc The Three Cases.}  Recall that non-line geodesics in $J^k$ come in three ``flavors'': heteroclinic, homoclinic and $x$-periodic. 
It is worth going into details around the interval $I$, the domain of the geodesic,   for each of the three cases.

($x$-Periodic). Choose   time origin so that $x(0) = x_0$ and $x(L/2) = x_1$.
Then $I = (0, L/2)$ or $(L/2, L)$ up to a period shift.  The minimizing
arcs correspond to half periods  of the $x$-periodic geodesic.  The  domain $\Omega$ projects onto the
interior   $(a,b)$ of the Hill interval.

(Heteroclinic.) If $\gamma$ is   heteroclinic
then $I =  \R$ and   
$c: \R \to \Omega$ is globally minimizing within $\Omega$.
If one or both endpoints $x_0, x_1$ is a local maximum of $F(x)^2$ then $\Omega$ projects to an interval
$(\alpha, \beta)$ strictly bigger  than $(x_0,x_1)$  That the interval $(\alpha, \beta)$ of a slab region
in the heteroclinic case is typically bigger than the corresponding Hill interval is essential in   the proof of theorem \ref{main}. 

(Homoclinic).  In this case the $x$ curve bounces once off the non-critical endpoint of the Hill
interval.  Say this interval is $b$ and that we translate time so that $x(0) = x_1$.
Then $I$ is of the form  $(-\infty, 0)$ or $(0, \infty)$.
The Hamilton-Jacobi minimality argument does not allow us to include $t = 0$
within the domain of $\gamma$ as $\gamma(0)$ is outside the open slab.  As to the  domain $\Omega$, 
it will project to either an interval $(\alpha, b)$  bigger than $(x_0,x_1)$  or project onto $(x_0,x_1)$, depending on whether the critical endpoint $x_0$ is 
a local maximum of $F(x)^2$ or not.

\vskip .3cm

{\bf Remark.}  {\it The global minimality of $\gamma$  within $\Omega$ persists
in the heteroclinic   turn-back case.}

\section{Magnetic  cut times. Periods II.}

\subsection{Definitions of cut and Maxwell times.} 

\begin{definition}
Let $\gamma: \R \to X$ be a geodesic in a length space (eg. a \sR manifold)
parameterized by arclength. 
\begin{itemize}
\item{}   The  cut time of  $\gamma$ is
$$ t_{cut} (\gamma)  := \sup \{ t>0 : \; \gamma|_{[0,t]} \;\; \text{is length-minimizing}  \}.$$

\item{} A  positive  time $t= t_{MAX}$ is called a Maxwell time for $\gamma$
if there is a geodesic distinct from $\gamma$ which connects $\gamma(0)$ to $\gamma(t_{MAX})$ and whose length is $t_{MAX}$.
We then call $\gamma(t_{MAX})$  a Maxwell point along $\gamma$.  
\end{itemize} 
\end{definition}

The `Maxwell time'' terminology is popular in the Russian literature but uncommon in the English literature.
We use it here, inspired by \cite{Ardentov1, Ardentov2, Ardentov3}.  
It is well-known that in \sR and \Ri metric spaces, geodesics fail
 to minimize when extended beyond their smallest   Maxwell time $t_{MAX}$.  
 Thus,  
  $$t_{cut}(\gamma)  \leq \inf\{ t:  t  \text{ is a Maxwell time for } \gamma \} $$
See for example, \cite{CheegerEbin}, Lemma 5.2, chapter 5 for the \Ri case. 
\subsection{Cut time and $x$-period} 
A first simple yet  important  result is:  

\begin{lemma}[Maxwell point, reflection argument]\label{Gen-Mart-cur}
If $F(x)$ is an even  polynomial, then any  $\R^3_F$-geodesic on which crosses $x = 0$ twice fails to minimize.
\end{lemma}
\begin{proof}
 If $F(x)$ is a even  polynomial, then $R(x,y,z) = (-x,y,z)$ is an isometry of $\mathbb{R}^3_{F}$. Let $c(t)$ be a $\R^3_F$-geodesic that crosses the plane $x=0$ twice, one at $A$ and another at $B$. The $c_1 (t):= R(c(t))$ also crosses $x=0$ at $A$ and $B$. Thus $B$ is a Maxwell point to $A$ along $c(t)$ and so $c(t)$ cannot minimize past $B$. 
\end{proof}

This lemma  says   $t_{MAX} \leq \frac{L(a,b)}{2}$ for    curves associated to even polynomials 
when the curve stars at   $x=0$.   
(See equation (\ref{eq: period}) for the integral expression for $L(a,b)$ where $G(x) = a + b F(x)$.)
We extend the lemma to hold regardless of starting point.  
\begin{figure}%
    {{\includegraphics[width=2.5cm]{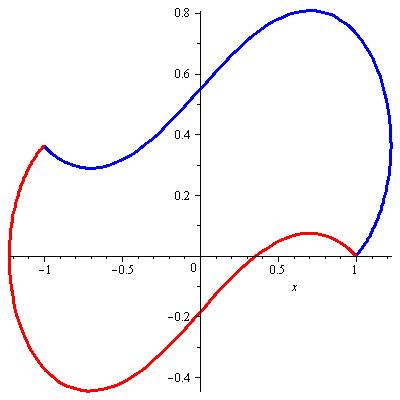} }}%
    \qquad
    {{\includegraphics[width=2.5cm]{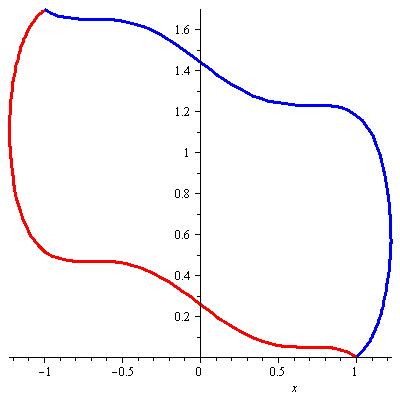} }}%
     \qquad
     {{\includegraphics[width=2.5cm]{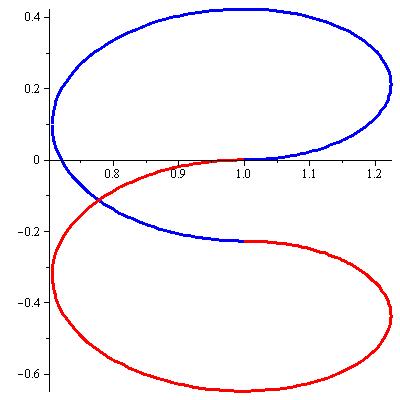} }}%
    \qquad
    {{\includegraphics[width=2.5cm]{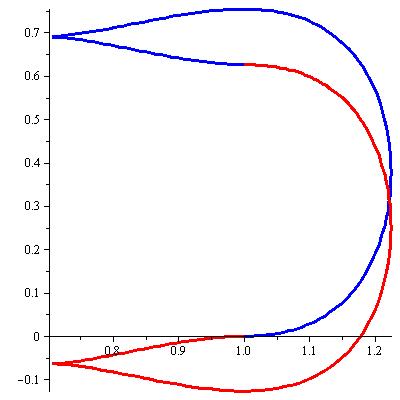} }}%
    \caption{The $(x-y)$-projections of  typical $x$-periodic geodesics indicating Maxwell points. 
    In each panel  
    two half-period curves share endpoints and are associated to the same
    polynomial $a + b F(x)$. } 
    \label{Equi-optimal} 
         \end{figure}

\begin{proposition}\label{Maxwell}
Let $c$ be a $x$-periodic geodesic on  $R^3 _F$
with $x$-period $L$. 
Then 

1.- $t_{cut}(c) \le L/2$ if  $F$ is even and $c$'s Hill interval   contains $0$.  

2.- $t_{cut}(c) \le L$ in all cases. 

\end{proposition}
\begin{proof}

 We start with the second case.  Let $c(t) = c_A (t) = (x_A (t), y_A (t), z_A (t)) $ be the geodesic,  
  let   $G(x)=a+bF(x)$ be its polynomial and   $[x_0,x_1]$ its Hill-interval.  Write $x_i= x_A (0)$. 
  If $x_i$ is interior to the Hill interval,  then there are exactly two magnetic geodesics passing
  through $c(0)$ and associated to $G(x)$, namely, the given one $c(t) = c_A(t)$
  and $c_B (t) = (x_B (t), y_B (t), z_B (t))$   characterized by  $\dot x_B (0) = - \dot x_A (0)$.  
  Then    $x_B (t) = x_A (-t)$ for all $t$. By $x$-periodicity we have $x_B (L) = x_A (L) = x_i$.   Proposition \ref{y-z as x function} tells us that 
   $c_A$ and $c_B$ have the same $y$ and $z$, periods,  $\Delta y, \Delta z$.  Thus 
$$ c_A(L) = c_A (0) + (0,  \Delta y, \Delta z) =  c_B(L).  $$ The geodesics  curves are distinct, showing that $L$ is a Maxwell time for $c$ and so $t_{cut} (c) \le L$.
In case $x_i$ is one of the Hill endpoints repeat the argument of Proof 1 of (i) of Theorem \ref{theorem-1} which was 
 given at the beginning of section \ref{sec: proof of B}  to conclude  that $c(L)$ is conjugate to $c(0)$ along $c$, so again $t_{cut}(c) \le L$.

We proceed now to  the first case where $F(x)$ is even.  Then  $G(x) = a + bF(x)$ is also even. 
By assumption the  Hill interval  for the $G$-geodesic $c(t)$ has the form $[-x_0,x_0]$ with $x_0 > 0$. 
Let us begin by assuming that $c(0) = (0, 0,0)$.  To determine the geodesic
we need the sign of $\dot x (0)$.  There are exactly two solutions, $(x(t), y(t), z(t))$ and
$(-x(t), y(t), z(t))$.  Both satisfy $c(-t) = -c(t)$ and in particular $x(-t) = - x(t)$.    Now if $L/2$ is the half-period of the $x$-curve we
have $x(L/2) = x(0) = 0$ but $\dot x (L/2) = - \dot x (0)$.  It follows that $x(t + L/2) = -x(t)$.  
One verifies that $y(t + L/2)$ and $y(t)$ both satisfy the differential equation $\dot y = G(x(t))$
from which it follows that $y(t + L/2) = y(t) + \Delta y$ with $\Delta y$ constant.
 Similarly $z (t + L/2) = z(t) + \Delta z$ with $\Delta z$ constant.
 
 {\it Attention!} The constants $\Delta y, \Delta z$ are exactly half the
 constants called $\Delta y, \Delta z$ in  proposition \ref{y-z as x function}.
 We can see this by writing out  $$y(t+ L) - y(t) = (y(t + L) - y(t + L/2)) + (y(t +L/2) - y(t)),$$
 and using the above half-period relation. An identical argument works for $\Delta z$.

 The general geodesic passing through $x =0$ at time $t = 0$
 has the form $c(t) + (0, \alpha, \beta)$ for $\alpha, \beta$ constants.
 Now any geodesic for $G$ is of the form $c(t + h)$ where $c(t)$ is as just described.
 It follows that every geodesic    for $G$  having Hill interval $[-x_0, x_0]$ satisfies
 the `monodromy relations'  \begin{equation}
 \label{half-period-relation}
 (x(t + L/2), y(t + L/2), z( t + L/2)) = (-x(t), y(t), z(t)) + (0, \Delta y , \Delta z).
 \end{equation}
 where
 \begin{equation}
 \Delta y  = \int_{-x_0} ^{x_0} \frac{G(x)}{\sqrt{ 1- G(x)^2}} dx, \qquad
 \Delta z = \int_{-x_0} ^{x_0} \frac{G(x)F(x)}{\sqrt{ 1- G(x)^2}} dx.
 \end{equation}
  Now, as described above,  there are exactly  two distinct $G$-geodesics passing through any point $(x_i, y_i, z_i)$ in $\R^3 _F$ provided $|x_i| < x_0$,
  namely  
 one heading right initially ($\dot x_i > 0$), and the other  heading left ($\dot x_i < 0$). By the above half-period identity,
  these two geodesics  re-intersect at the same point $(-x_i, y_i + \Delta y, z_i + \Delta z)$ a   time $L/2$ later.
  Consequently $L/2$ is a Maxwell time. 
\end{proof}

Figure  \ref{Equi-optimal} illustrates the Proposition by showing
the x-y projections of the two geodesics sharing endpoints
for several  polynomials $G(x) = a + b F(x)$.  

We make a conjecture concerning a property that  \cite{Ardentov1, Ardentov2, Ardentov3}
call ``equi-optimality'' which proved useful both technically and organizationally for their proofs.
 
\begin{definition}  We say that the arc-length parameterized geodesic
$\gamma: \R \to X$ is {\it equi-optimal} if its cut-lengths are independent of where
we start on the geodesic. In other words, for any real $s$, let   $\gamma_{s} (t) = \gamma (t-s)$ be  the time translated version of $\gamma$, having new starting point $\gamma_s (0) = \gamma(s)$. 
Then $\gamma$ is   equi-optimal if   $t_{cut}(\gamma_{s})$ is independent of $s$. 

We say that a length space is equi-optimal if all the geodesics are equi-optimal.
\end{definition}

{\sc Conjecture.}  $J^k$ and $\R^{3}_F$ are equi-optimal.

This conjecture is well-known to hold for $J^1$,  the Heisenberg group.
For $J^2$,   the Engel group, the conjecture was established in \cite{Sahkov} by   computations with elliptic functions.  
The work presented here suggests the conjecture might hold for all $J^k$ but we are far from a proof.

This proposition almost proves the conjecture on equi-optimality
for geodesics in the case that  $F(x)$ is even and its  Hill interval is $[-x_0, x_0]$, for some $x_0 > 0$.   
Missing is a proof that $t_{cut}(c) = L/2$: that is,  that  no Maxwell or conjugate times can be less than $L/2$.
If the starting point is one of the Hill endpoints $\pm x_0$ then this equality for $t_{cut}$ follows by the Hamilton-Jacobi argument, proposition \ref{HJminimality}.
However, we do not know how to get this missing piece to the proof when $x(0)$ is  interior to the Hill interval.

  \vskip .3cm
  
  \subsection{Periods II} 
  \label{subsec: Periods} 
Under the assumption that $F$ is even
  and our Hill interval contains $0$,  we set
  \begin{eqnarray}
 \label{half-periods-2}
\Delta y (a,b) & = & 2 \int_0 ^{u(a,b)}   \frac{G(x)}{\sqrt{ 1- G(x)^2}} dx, \\
 \Delta z  (a,b) &= &2 \int_0 ^{u(a,b)} \frac{G(x)F(x)}{\sqrt{ 1- G(x)^2}} dx,
 \end{eqnarray}
  where   $G(x)   = a + b F(x)$  
 and   $u  = u (a,b)$ is the first positive solution to $ G(x)^2 = 1$ 

These are the translations suffered by $y(t), z(t)$ after each half period.
Compare equation (\ref{half-period-relation}), for the half-period recall that
\begin{equation}
\Delta t (a,b) = 2  \int_0 ^{u(a,b)} \frac{dx}{\sqrt{ 1- G(x)^2}} . 
\end{equation}  This half-period is the {\it length} of the geodesic over a
half-period and equals $L(a,b)/2$ where $L(a,b)$ is the period, observe that
$$ |\Delta y(a,b)|  < \Delta t (a,b),$$
since $|G(x)| < 1$ on $(0,u)$.

 {\bf We collectively refer to $\Delta t (a,b),  \Delta y (a,b), \Delta z (a,b)$
as the periods associated to $(a,b)$}. It will be crucial that 
they are independent of where we start  
along the curve, i.e. of the   $t$ in equation  (\ref{half-period-relation}).
 The  functions $\Delta y (a,b), \Delta z (a,b)$ and $\Delta t(a,b) = L(a,b)/2$ are analytic functions of $(a,b)$
in a neighborhood of any value   $(a,b)$ for which they are finite.

The  periods at $(a,b)$ are finite if   the Hill endpoint  $u = u (a, b)$ is  a simple root of $1 - G(x)^2$. ( Note that,
by definition of ``Hill interval'',   $1 - G(x)^2$ has no zeros in the interior  $(-u, u)$ of its Hill interval.)  
 The endpoint $u$ is a  double root if and only if it is a critical point of $G$ in which case
 the geodesic is not periodic.     Situations where new roots appear  in the interior of the Hill interval and where  $u \to \infty$ arise through limits which appear when we investigate candidate
 long-period minimizers approaching a heteroclinic geodesic  at a   key step below  in our argument for proving theorem \ref{main}.

 \section{Proof of  theorem \ref{main} } \label{sec:proof-main-the}
 
In this  long technical section  we will prove
theorem \ref{main}.   We  begin  with an   outline
 for the   section and hence for  the   proof.

\subsection{Outline} Theorem  \ref{main} asserts that the heteroclinic geodesic  for a certain class of seagull potentials $F$,
when projected to the magnetic space for $F$,  is a global minimizer 
there, and hence the   geodesic itself is a  global minimizer on $J^k$.  
By the `base geodesic' we will mean  this  projected geodesic to the magnetic space for $F$.
To show that  the base geodesic globally minimizes in the magnetic space we proceed   by contradiction.  All geodesics in the magnetic space
are governed by   polynomials $G(x)  = a + b F(x)$  in the pencil for $F$. If the base geodesic is
not globally minimal then we can find a sequence of endpoints symmetrically placed 
along the base geodesic whose distance from each other  tends to infinity and a sequence of  $G$'s whose geodesic  arcs share
these endpoints and which are shorter,  or at least no longer, than the corresponding arc of the base geodesic.  
 An application of the Hamilton-Jacobi method shows that  these shorter  geodesics must leave the
slab $-\beta \le x \le \beta$ which strictly contains the slab $-1 \le x \le 1$ of  the base geodesic.   These two pieces of information - the endpoint conditions combined with the leaving-of-the-slab 
 yield a compactness
for the polynomial family $G$:  namely we must have that $|a + b| \le 1$ and $|a-b| \le 1$. We  call this locus of points in the $a,b$ plane   the Diamond, denoted by $DIAM$ below.   See definition \ref{competing} which acts as a kind of summary.   All  this is done in the next `set-up' subsection \ref{subsec:set-up}.

Subsection \ref{subsec: period analysis 2}  continues the analysis of  endpoints
and period asymptotics begun in the previous section  and introduces one of our key tools, the y and z ``costs''.
See  equations
(\ref{costy1}).   These are the  differences of two periods associated to the competing geodesics coming from the Diamond,
and can be thought of as `renormalized' periods.  We note    that  as the endpoints tend to infinity on our base geodesic, the coefficients $a,b$ encoding the competing geodesics
must   tend to a ``Z'' - the union of three line segments -contained in the Diamond.  We denote
these segments by `Leg 1', `Leg 2' and `Leg 3'.   See figure \ref{fig:PEN-picture}. 
In this way we reduce the work to that of  understanding the asymptotics of  the y and z costs  as we approach these three legs.  We end the subsection with  
the statement  of  Proposition \ref{prop: period asymptotics}, 
a proposition on this asymptotics which   almost immediately yields  the Theorem.  
 
Subsections \ref{subsec:geting-rid} and \ref{subsec:get-rid-2}  are devoted to eliminating the three Legs of the Z  one at a time -  thus showing that
 the competing  geodesics cannot be simultaneously shorter than the base geodesic  and share its endpoints.   
The method of elimination is   essentially calculus, through the computing  the asymptotics and the variations of the y and costs.  In the short final subsection \ref{subsec: end proof main}
we  show how Proposition \ref{prop: period asymptotics} implies the Theorem.

\vskip .3cm
 
\subsection{Set up for the proof}
\label{subsec:set-up} 
 
Recall (definition \ref{def: seagull})  that a seagull potential is even, achieves its
global maximum value of   $1$ at $x = \pm a$ and satisfies $0 < F(x) < 1$
for $-a < x < a$.  Moreover $F'(x) < 0$ for $x > a$ so that  $F$  tends to $-\infty$ as $x \to \infty$. 
It follows that $F$'s Hill intervals are $[-\beta, a], [-a, a], [a,\beta]$ where
$x = \beta$ is the unique positive $x$ having $F(x) = -1$. 
Use a scaling symmetry $x \mapsto hx$ to scale our seagull potential  $F$
in order to place the maximum points $x = \pm a$ at   $\pm 1$.   

Our claimed globally minimizing geodesic -  is the  direct heteroclinic geodesic $\gamma_0: \R \to J^k$ 
for $F$ with Hill  interval   $[-1,1]$ and whose $x$ curve is monotonic increasing.
Thus its   $x$-curve limits to $-1$ in backward time and to $x=+1$ in forward time.
 Let   $c_0 = \pi_F \circ \gamma_0: \R \to \R^3 _F$ be  its projection to the plane -the curve referred to as
 the `base geodesic' above. 
  It suffices to show that   $c_0$ is a globally minimizing in $\R^3 _F$ to conclude Theorem \ref{main}. 
 
Next, we assume that $\beta = \sqrt 3$, which is to say  $F(\sqrt{3}) = -1$.
 It follows   that the other Hill intervals for
 $F$ are $[1, \sqrt{3}]$ and $[-\sqrt{3}, -1]$. 
   It follows from Proposition \ref{HJ_min} and the discussion around it that we  have a smooth solution to 
 the Hamilton-Jacobi equation for $F$  on the slab domain 
 $$\Omega:= \{ (x, y, z):  -\sqrt{3} < x < \sqrt{3}  \}.$$ 
It   follows  from  Proposition \ref{HJ_min}   that $c_0$ globally minimizes
 within $\Omega$.

We now argue by contradiction. If $c_0$ fails to globally minimize, then  there   exist a family of   shorter geodesics in $\R^3 _F$   connecting distant  endpoints of  $c_0$. Due to the Hamilton-Jacobi result just discussed, these   shorter geodesics   must all leave $\Omega$.    
 The proof will be  completed by   showing that these shorter geodesics cannot exist.

To begin, we argue that these shorter geodesics must be arcs of periodic geodesics.  To this purpose, and
to simplify book-keeping, we  shift the time origin and translate as needed so that $c(0) = (0,0,0)$.
Write $$c_0 (t) = (x_0 (t), y_0 (t), z_0 (t)).$$ 
It follows from the evenness of $F(x)$ that  
$c_0 (-t) = (-x_0 (t), -y_0 (t) , - z_0 (t))$.  Since  $c_0$ fails to globally  minimize, we have, for all $T/2$ sufficiently large,
a shorter geodesic $c= c_T$  in $\R^3 _F$ joining $c_0 (-T/2)$ to $c_0 (T/2)$.   

Set
$$\delta = x_0 (T/2),$$
so that for $T$ large $\delta$ is very close to $1$.  
$T$ and $\delta$ are related,
through the $x$-differential equation, by:  
 \beq
 \label{eq: delta}
 T  = 2  \int_{0} ^{\delta} \frac{dx}{\sqrt{ 1- F(x)^2}} .
 \eeq
 
\begin{figure}%
     {{\includegraphics[width=3cm]{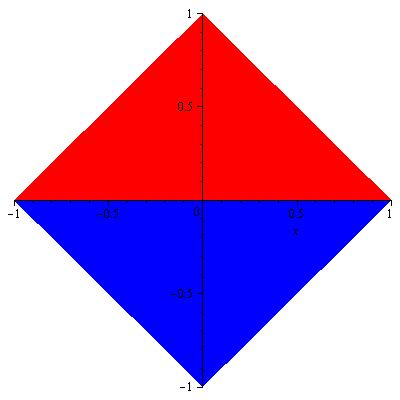} }}%
      {{\includegraphics[width=3cm]{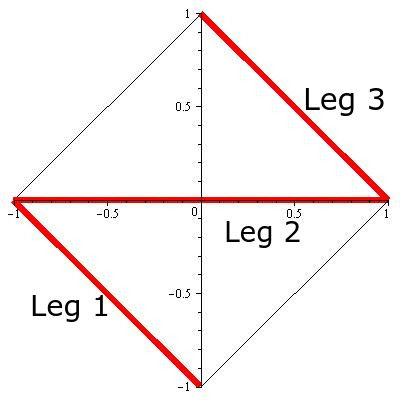} }}%
      \caption{ On the left panel is the  diamond, DIAM, whose points parameterize competing geodesics.
      On the `Z' along whose 3 line segments, denoted $Leg 1, Leg 2$, and $Leg 3$,  one or more of the periods blow up.  The 
      coordinates of the axes are  $a$ and $b$ of  ``$G(x) = a + b F(x)$''.}
   \label{fig:PEN-picture}   
   \end{figure}

The allegedly shorter geodesics $c = c_T = (x(t), y(t), z(t))$, being a geodesic on $\R^3 _F$ is associated
to some polynomial  $ G (x)  = a + b F (x)$ in the pencil of $F$.     
Since endpoints of $c$ and $c_0 ([-T/2, T/2])$  match up, the geodesic $c$
starts at $x = -\delta = x_0 (-T/2)$ close to $-1$, crosses $x= 0$ and ends up at $x = \delta$ close to $+1$,
while in doing so it  must leave   the   slab.
Thus its $x(t)$  must reach a maximum point $u \ge \sqrt{3}$, or minimum point $u < -\sqrt{3}$ and 
return to $x_0 (T/2) = \delta < 1$.     It follows that $x(t)$ when extended to $t \in \R$ is  periodic,
with Hill interval $[-u,u]$, $u \ge \sqrt{3}$. In more detail: since $x(t)$ returns from $u$ we   have
$G'(u) \ne 0$.  By evenness of $G(x)$ and the fact that $x(t)$ crosses zero, the Hill interval associated to $c$   is   $[-u,u]$.
 In particular, since both  $1$ and $\sqrt{3}$ must
be in the   Hill interval of $x(t)$ and since   $F(1) = 1,  F(\sqrt{3}) = -1$, by evaluating $G$ at
these two points  we have that
 $|a+b|  < 1$ and 
$|a-b| \le 1$.  The latter is an inequality since $u = \sqrt{3}$ is allowed, while the former is
strict since $1$ must be in the interior of the Hill interval of $G$ in order for the $x$-curve to make it all
the way to $\sqrt{3}$. We will denote the region that we just describe in the following way, 
 \beq
 DIAM: = \{  (a,b):  |a + b| < 1 \text{ and } |a-b| \le  1, b \ne 0 \}.
\label{eq: diamond} 
\eeq  
DIAM is short for ``diamond''. We call $DIAM_+$ and $DIAM_{-}$ to the point $(a,b)$ in $DIAM$ such that $0<b$ and $b<0$, respectively.  See figure \ref{fig:PEN-picture}.

 We just argued that we need only compare $c_0$ to those geodesics arising from  
 $G = a+ b F(x)$ with  $(a,b)$ lying in $DIAM$ in the $a-b$ plane.
Let us formalize the above discusion by: 
\begin{definition}
\label{competing}  [The family of competing geodesics] For each  point 
$(a,b)$ in the square $DIAM$ and each $T >0 $   let $c = c_{a,b;T}: [0, L(a,b)/2] \to \R^3 _F$ be the   geodesic arc  for $G(x) = a + b F(x)$ which starts at $t = 0$ at  the point $c(0)=  c_0 (-T/2) $
so that $x(0) = -\delta$,  
and has  $\dot x >0$, followed  for half of its $x$-period,   $\Delta t (a,b) = L(a,b)/2$.
In this way its $x$ curve achieves a maximum of $u = u(a,b) \ge \sqrt{3}$
and ends when   $x =\delta = +x_0 (T/2)$  at time   
$$\Delta t (a,b) := L(a,b)/2.$$    
\end{definition} 
   
   {\bf Remarks.} 1.  We do not have to worry about $x$-periodic curves which hit both $u$ and $-u$
   since these are longer than their half-period $L/2$ and hence do not minimize by (i) of Proposition 
   \ref{Maxwell}. 
   
   2.  If, in the definition, we took $\dot x < 0$ at time $t =0$ instead, we would end up with the other
   $G$-geodesic, again for a half-period $L/2$.  See again  the arguments around (i) of
   Proposition \ref{Maxwell} or figure \ref{Equi-optimal}.  This other curve has the same endpoints and same
   length as our curve, so is identical for all our purposes.  The choice $\dot x > 0$ is just made 
   so as to simplify the exposition. 

\vskip .3cm 

\subsection{Period asymptotics, continued}
\label{subsec: period analysis 2}
For    our   competing geodesics $c$, as described in definition \ref{competing},  the $x$-values of both endpoints match the $x$-values   of $c_0$
by design.   We will complete the proof of the Theorem  by  showing that  if
their  y and z values match  $c_0$'s 
then $c$ is longer: $T < L(a,b)/2$. 

To this purpose  recall the half-period relations for such periodic geodesics,  
equations (\ref{half-period-relation})  and (\ref{half-periods-2}) for how $y$ and $z$ change in a half-period, these assert that
$$c(\Delta t (a,b)) = ( + \delta,   -y_0 (T/2),  - z_0 (T/2)) + (0, \Delta y(a,b), \Delta z (a,b).$$
The requirement that the far  endpoint of $c = c_{a,b,T}$ agrees with endpoint of $c_0 ([-T/2, T/2])$ 
is  the  requirement that their y and z periods $\Delta y (a,b), \Delta z (a,b)$ as given
by (\ref{half-periods-2}) satisfy  
 \begin{eqnarray}
 \Delta y(a,b) &=& \Delta y_0 (T), \\
  \Delta z(a,b) &=& \Delta z_0 (T), 
   \label{endpoints} 
 \end{eqnarray}
 where 
 \begin{eqnarray*}
 \Delta y_0 (T)    &=& 2  \int_{0} ^{\delta} \frac{F(x)}{\sqrt{ 1- F(x)^2}} dx, \\
 \Delta z_0 (T)   &= & 2 \int_{0} ^{\delta} \frac{F(x)^2}{\sqrt{ 1- F(x)^2}} dx,
 \end{eqnarray*}
 are the   corresponding changes in $y$ and $z$ suffered as we travel our heteroclinic orbit
 from  $c_0 (-T/2)$ to $c_0 (T/2)$,
 namely   $\Delta y_0 (T) = 2y_0 (T/2)$ and
 $\Delta z_0 (T) =2 z_0 (T/2)$).

 {\bf Remark.}  $\Delta y (a,b)$ and $\Delta z (a,b)$
 are independent of the starting point on the curve $c$.  Compare equation (\ref{half-period-relation}).
 This $T$-independence simplifies arguments in an   
 essential way and is   where we use the assumption that  $F$ is even.
 
As $T \to \infty$ we have  $\delta (T) \to 1$ and 
$\Delta y_0 (T), \Delta z_0 (T) \to + \infty$.   
Thus, in requiring  the endpoint conditions, equations (\ref{endpoints}), to hold for our competing curves
we are forced to investigate the periods $\Delta y(a,b), \Delta z(a,b)$ as they tend to infinity.
 All three   periods $\Delta t, \Delta y, \Delta z$ are analytic  functions of $(a,b) \in DIAM$  away from where they blow up.     Analysis of the (a,b) periods  near the blow-up loci 
of equations (\ref{half-period-relation})  and (\ref{half-periods-2}) is the whole game now.  

These periods blow up   along  three line segments in the diamond and nowhere else.
These segments    form  a    tilted ``Z''  whose middle segment,
labelled $Leg 2$,  is the  segment of the $a$-axis inside the   diamond,  see figure \ref{fig:PEN-picture}.  
 The other two strokes of the $Z$  are made up of the bounding edges of the diamond, denoted by $Leg 1$, the locus
 $a+ b = -1$  and $Leg 3$, the locus $a+ b = +1$.   Blow up of periods  along $Leg 2$,  
($b =0$),  occurs since the 
Hill endpoint   $u(a,b)  \to \infty$ as $b \to 0$.
Blow up along $Leg 1$ and $Leg 3$ occur because as we approach either of these legs 
 we have that  $G(1)^2 \to 1$, that is, the   `mountain peaks'' $x = \pm 1$ of $G(x)^2$  get closer and closer to satisfying  $G(x)^2 = 1$, finally
   touching  $G(x)^2 = 1$ at points of these legs,  making periods
very long  through the reciprocal $1/\sqrt{1 - G(x)^2}$   occuring  in all three period integrals.   

{\bf Remark.} We cannot get 
  the other critical point $x = 0$ of $G(x)$ to attain the level $G(x)^2 = 1$  
  in the closure of the diamond and thereby lead to divergent periods,  except possibly at the points $(\pm 1, 0)$
  which will be dealt with directly.    Indeed  since $|F(0)| < 1$
the condition $G(0)^2 = 1$ which is  $(a + b F(0))^2 = 1$ together with $|a \pm b| \le 1$ yields $a = \pm 1$.

 We focus on the differences in periods
\begin{eqnarray}
Cost_y (a,b) &=& \Delta t (a,b) - \Delta y (a,b),\\
Cost_z (a,b)  &=& \Delta t (a,b) - \Delta z (a,b),
\label{costy1} 
\end{eqnarray} 
rather than  the periods themselves.  The advantage gained is that these difference of   periods
  have finite limits as we tend to $Leg 1$ and $Leg 3$ (except for $(\pm 1, 0)$)  and so   extend to continuous functions on the closure of the entire Diamond minus $Leg 2$ i.e. minus the $a$-axis.   
 We will compare these `renormalized' periods to the analogous   quantities
for the heteroclinic geodesic : 
\begin{eqnarray*}
Cost_{0,y}  (T) & = & T - \Delta y_0 (T)  \\
   & = &  2 \int_0 ^{\delta}  \frac{ (1 - F(x))}{\sqrt{ 1 - F(x)^2 }} dx,
   \end{eqnarray*} 
   and
   \begin{eqnarray*}
Cost_{0,z}  (T) & = & T - \Delta z_0 (T)  \\
   & = &  2 \int_0 ^{\delta}  \frac{ (1 - F^2(x))}{\sqrt{ 1 - F(x)^2 }} dx,
   \end{eqnarray*} 
   which both have finite limits as $T \to \infty$, that is to say as $\delta \to 1$. 
   
   \begin{lemma}  \label{lem: cost at infinity}The functions
   $Cost_{0,y}  (T), Cost_{0,z}  (T)$ are strictly monotone increasing
   in $T$ and tend to finite positive limits as $T \to \infty$:  
    $$Cost_{0, y} (\infty) = \lim_{T \to \infty}  Cost_{0,y}  (T),$$
   and 
   $$Cost_{0,z}(\infty) =  \lim_{T \to \infty}  Cost_{0,z}  (T).$$ 
   These limits  can be obtained by setting $\delta =1$ in the integral expressions given just above for $Cost_{0,y}  (T)$ and $Cost_{0,z}  (T)$  
  \end{lemma}  
   
  { Proof.} The integrands are all positive and behave like 
   $ \sqrt{|1-x|}$ near $x =1$.

   The proof of Theorem \ref{main} will be completed upon establishing the following result: 
   \begin{proposition}
   \label{prop: period asymptotics} 
   
   a) In $DIAM_+$, including all along $Leg 3$  we have   
   $$Cost_y (a,b) > Cost_{0, y} (\infty).$$
   
   b) $Cost_z (a,b) \to + \infty$ as  we approach any point of either   $Leg 1$ or $Leg 2$
   along curves in  $DIAM_-$, including the  point $(1,0)$.    
   \end{proposition}
  
\subsection{Getting rid  $Leg 1$ and $Leg 2$}\label{subsec:geting-rid}
      
 We will  prove  proposition \ref{prop: period asymptotics} by obtaining   detailed asymptotics for the costs (periods)  close to the Z. We will split the proof into two main parts, according to the two parts of the proposition.  Part (b) for points in    $DIAM_-$ 
itself splits into  three cases, labelled below as  ``getting rid of points near $Leg 1$", 
``getting rid of points near $(a,b) = (1,0)$ with $b < 0$" and "getting 
rid of points near $Leg 2$  with $b < 0$". Part (a) for $DIAM_+$ is presented as a single 
case ``getting rid of points on $Leg 3$" whose proof consists of    two steps.

  
\subsubsection{ Getting rid of points near $Leg 1$}

On and near  $Leg 1$ we have   $a, b < 0$.  Since $b < 0$ the absolute minimum of $G = a + b F(x)$
occurs when $x=1$ and its value there is  $a + b$ which is $-1$ at points of $Leg 1$.  
The integral computing $\Delta y (a,b)$  is that of  $G(x)/ \sqrt{ 1- G(x)^2} $ over $[0, u]$
where   $u = u(a, b)$ is the first positive solution to $G(x) = 1$. The denominator 
$ \sqrt{ 1- G(x)^2} $ goes to zero at $x =1$ and at $x = u$.  Its behaviour near $x =u$ is like $1/\sqrt{|u-x|}$
which is integrable.  Its behaviour near $x = 1$ is like $-1/|1-x|$
so the integral  l diverges logarithmically to $-\infty$.
It follows that as $(a,b)$ tends to any point of $Leg 1$. 
we have that $\Delta y(a,b) \to - \infty$. 
It follows that for all $(a,b)$ sufficiently close to $Leg 1$ we have $\Delta y < 0$.
On the other hand,  $\Delta y_0 (T) \to + \infty$ with $T$,  the divergence 
being due to the behaviour of   the integrand
near $x=1$.  This makes it  impossible to satisfy 
the endpoint conditions (\ref{endpoints}.)  

We have not yet   proved   claim (a) of the proposition. Note that $\Delta t > 0$,
so by our just established asymptotics for $\Delta y$ we have that $Cost_y \to + \infty$ as we approach $Leg 1$, while
in comparisont $Cost_{0,y} (\infty)$ is finite.   

Regarding the claimed behaviour of $Cost_z (a,b)$
near $Leg 1$ in the proposition,  use that  $F(1) = +1$ so that an identical analysis applied to the integral
expression for $\Delta z(a,b)$ shows that
$\Delta z (a,b) \to -\infty$ as $(a,b)$ tends to any point of $Leg 1$.  On the other hand
$\Delta t (a,b) > 0$  so that again
$Cost_z (a,b) = \Delta t (a,b) - \Delta z (a,b) \to + \infty$.   

\vskip .3cm 

\subsubsection{ Getting rid of  points near $(1,0)$ with $b < 0$} 

We begin by  parameterizing  the lower diamond by using lines along which the
upper bound $u$ of integration is constant.  Since $b < 0$ we have that $G$ rises from its global minimum
of $a + b <1 $ which occurs at  $x = 1$,  up past the value $a-b$ at $x = \sqrt{3}$, and on  until it hits the unique    positive
$u$ such that $G(u) = 1$.  Set $a = 1- \tau$. Then the equation $G(u) =1$ says 
$(1 - \tau) + b F(u) = 1$ which is to say $b = \tau/F(u)$.  We then parameterize the points of $DIAM_-$
by
$$a = 1 - \tau,$$
$$b = \tau/F(u), u > \sqrt{3}.$$
In this parameterization, by  fixing $u$ and varying $\tau>0 $ down to $0$ we
approach $(1,0)$ along lines through $(1,0)$ having slope $1/|F(u)|$.  Note as $u$
increases from $\sqrt{3}$ to $u = \infty$ the slope of these lines goes from $1$ to $0$,
thus sweeping out all of $DIAM_-$.  

In these coordinates
$$G = (1 - \tau) + \frac{\tau}{F(u)} F(x),$$
we compute that
$$1 - FG = 1 - F + \tau F -  \frac{\tau}{F(u)} F(x)^2, $$
which tends to $1-F$ as $\tau \to 0$.

On the other hand,  $1 - G^2 = (1-G)(1+ G)$
and  $1+G \le 2$ on $[0, u]$
so that
$$(1 - G^2) \le 2(1-G),$$
on this interval.  Now observe that for $x \in [0, u]$ we have 
\begin{equation}
\begin{split}
 1-G & =  \tau -  \frac{\tau}{F(u)} F(x) \\
 & =   \frac{\tau}{|F(u)|}  ( F(x) - F(u)),  
\end{split}
\end{equation}  
where we have used that $F(u) <  0$. Note that this expression is   positive  in $[0,u]$ since  $F(x) > F(u)$ for $0 < x < u$.
Thus  for $x \in [0, u]$ we have 
$$(1-G(x)^2) \le 2 \frac{\tau}{|F(u)|}  ( F(x) - F(u)), $$
 and hence that  
$$\frac{1}{\sqrt{1 - G(x)^2}} \ge \frac{|F(u)|^{1/2}}{\sqrt{2 \tau}} \frac{1}{\sqrt{F(x) - F(u)}} .$$
Our integrand for $Cost_z (a,b)$ is $\frac{1- FG} {\sqrt{1 - G(x)^2}}$. 
Since $(1-F(x)) \ge 0$ everywhere, our expansion of $(1-FG)$ above 
yields   
$$\frac{1 - F(x) G(x)}{\sqrt{1 - G(x)^2}} \ge \{ (1 - F)  - \tau |F(x)| + \frac{\tau}{|F(u)|} F(x)^2 \} \frac{1}{\sqrt{2 \tau}}
 \frac{|F(u)|^{1/2}}{\sqrt{F(x) - F(u)}}. $$

The three terms on the right hand side that 
 we get by freezing $\tau$ and $u$, namely $(1- F(x))/{\sqrt{F(x) - F(u)}}$,  $|F(x)|/{\sqrt{F(x) - F(u)}}$
and  $F(x)^2/{\sqrt{F(x) - F(u)}}$ all have finite integrals over $[0,u]$.   After integration,  let 
  $\tau$ vary down to zero, for   frozen $u$.  The last two terms go to zero like $\sqrt{\tau}$.
The first term goes to $+ \infty$ like $1/\sqrt{\tau}$.  

We have shown that approaching $(1,0)$ along any line in $DIAM_-$ the value of $Cost_z (a,b)$
approaches $+ \infty$.

\vskip .3cm

\subsubsection{ Getting rid of  points near $Leg 2$  with $b < 0$} 

We will show that $\Delta z(a, -\eps) \to -\infty$ as $\eps \to 0$ with $\eps >0$, and   for all $a$ with 
$-1 < a < 1$.

We start with the case $a> 0, b = - \eps$ with $\eps >0$ going to zero.
Since $F \le 1  $ we have that
$$G(x)  = a - \eps F(x)  \ge a- \eps,$$ 
and in particular, for $\eps$ small enough $G>0$.
It follows that   $G(x) F(x)$ has the same sign  as $F(x)$.  
 Let $z_0$ be the first positive zero of $F(x)$ so that $ G F \ge 0$ on $[0, z_0]$ 
while $G(x)F(x) < 0$ on $(z_0, u]$ where $u$ is the first positive solution to $G(x) = +1$. 
(Note $1 < z_0 < \sqrt{3}$.) 
 Then  
\begin{equation}
\begin{split}
\Delta z (a, -\eps) &= 2\int_0 ^{z_0} \frac{G(x) F(x)}{ \sqrt{ 1- G(x)^2}} dx + 2\int_{z_0} ^{u} \frac{G(x) F(x)}{ \sqrt{ 1- G(x)^2}} dx \\
& =  I_0 (\eps) + I_1 (\eps).
\end{split}
\end{equation}
The first integral $I_0 (\eps)$ is positive and the second is negative. 
As $\eps \to 0$ the first integral tends to a finite     positive value. 
Indeed on $[0, z_0]$ we have $a- \eps \le G(x) \le a$ so that as $\eps \to 0$ we have
$GF/ \sqrt{ 1- G(x)} \to a F/ \sqrt{ 1 - a^2} $ which is a (bounded!) polynomial on $[0, z_0]$,
leading to a finite limit for $I_0 (\eps)$.  

Now
$$ I_1 (\eps) =- 2 \int_{z_0} ^u \frac{ G(x) |F(x)|} {\sqrt{ 1 - G(x)^2}} dx,$$
and $G(x)$ is monotonically increasing from $a$ to
$1$ on $[z_0 , u]$ so that $G \ge a$ while $1/\sqrt{1-G^2} \ge 1/ \sqrt{ 1- a^2}$,
so 
$$\frac{G(x)}{\sqrt{ 1 - G(x)^2}} \ge \frac{a}{ \sqrt{1 - a^2}} \;\;  \text{and}\;\;
   \frac{ G(x) |F(x)|}{\sqrt{ 1 - G(x)^2}}  \ge \frac{a|F(x)|}{ \sqrt{1 - a^2}} $$
 on this interval. 
Now we will need estimate (\ref{eq: root asymptotics}) below for the   Hill endpoint $u = u(a, -\eps)$.
We obtain the estimate  by approximately  solving
$a - \eps F(x) = 1$ or $$-F(u) = (1-a)/\eps.$$
  Since $\eps$ is very small we are solving
a polynomial equation   $p(x) = w$ for $w = (1-a)/\eps >> 1$, where our polynomial $p(x)$
is $-F(x)$ which has the form 
 $p(x) = a_0 x^{2k} + \ldots $
and in particular has   degree $2k$ with $a_0 > 0$.  The solution can be expanded as
\begin{equation}
u(a, -\eps) = (\frac{1-a}{a_0})^{1/2k} \frac{1}{\eps^{1/2k}} +  c + O(\eps^{1/2k}) ,
\label{eq: root asymptotics}
\end{equation}
valid as $\eps \to 0$.  Here   $c$ is constant.  
 See Appendix B for this standard result regarding asymptotically solving real polynomials.  
Moreover we have that     $|F(x)| > A x^{2k}$ for some $A > 0$
(e.g. $A = a_0 - \eps$) for all $x$ sufficiently large.
Thus 
$$I_1 \le - \frac{2a}{\sqrt{1 - a^2}}  \int_{z_0} ^{u(a, -\eps)} A x^{2k} dx.$$
The last integral yields $-C u^{2k+1} +O(1) = -\tilde C  \frac{1}{ \eps^{1 + 1/2k}}+ O(1) $
with $C, \tilde C$ positive constants, showing that   $I_1 (\eps) \to -\infty$ as $\eps \to 0$.

For the case $a \leq  0, b = -\eps$, we define $u_0 = u_0 (\eps)$ and $u_1 = u_1 (\eps)$ as the positive numbers where $G(u_0) = 0$ and $G(u_1) = 1$ then $G(x)$ is negative on $[0,u_0)$ and positive on $(u_0,u_1]$. We can take $\epsilon$ small enough to have $z_0 < u_0$, where $z_0$ is again the positive number such that $F(z_0) =0$, then $F(x)G(x)$ is negative on $[0,z_0]$ and $[u_0,u_1]$, while, $F(x)G(x)$ is positive on $[z_0,u_0]$. We  use  the same method but  replace the point $z_0$ by $u_0$, that is,  we split the integral  into $I_0(\epsilon)$ and $I_1(\epsilon)$, where $I_0$ is over $[0,u_0]$ and $I_1$   over  $[u_0,u_1]$. For small enough $\epsilon$ we have that $I_0(\epsilon)$ is positive and $I_1(\epsilon)$ is   negative. 
  
We estimate  $u_1$ and $u_2$ as before  by approximately solving $a - \eps F(x) = 0$ and $a - \eps F(x) = 1$. Again setting   $-F(x)$ to  $p(x) = a_0 x^{2k} + \ldots $with $a_0 > 0$ we find     
$$ u_0(\eps) = (\frac{-a}{a_0})^{1/2k} \frac{1}{\eps^{1/2k}}+ c_0 + O(\eps^{1/2k})\;\; \text{and}\;\; u_1(\eps) = (\frac{1-a}{a_0})^{1/2k} \frac{1}{\eps^{1/2k}}+ c_1 + O(\eps^{1/2k}).$$ 
Here $c_0$ and $c_1$ are constant. We have that    $I_0(\epsilon)$ and $I_1(\epsilon)$  now  go to infinity when $\eps$ goes to $0$ but because $\frac{-a}{a_0} < \frac{1-a}{a_0}$, we have that $I_1$  dominates $I_0$
so their sum again goes to negative infinity.

 \vskip .3cm 
 
\subsection{Getting rid of  $Leg 3$} \label{subsec:get-rid-2}
 
At this point we need to be precise regarding the class of polynomials.
Here is the promised definition.

\begin{definition} 
\label{def: specific class} The specific class of polynomials for the theorem
are obtained by taking $a =1$  and   $W(x)$ of the form 
$$W(x) = \frac{1}{6}(1+P(x)+(\frac{x}{\sqrt{3}})^{14})$$
in the equation (\ref{eq: W}) defining $F(x)$.
We insist that   $P(x)$ has degree at most $14$,   $P(\sqrt{3}) = 1$ 
and  
\begin{itemize} 
\item{i)} $ (\frac{x}{\sqrt{3}})^{14} \leq |P(x)| \leq 1$ if $|x| \leq \sqrt{3}$,
\item{ii)} $1\leq P(x) \leq (\frac{x}{\sqrt{3}})^{14}$ if $\sqrt{3} \leq |x|$,
\item{iii)} $|P'(x)| < \frac{14}{\sqrt{3}}| (\frac{x}{\sqrt{3}})^{13}|$ if $\sqrt{3} \leq |x|$.
\end{itemize} 
\end{definition} 
All such $F$'s have  Hill region $[-\sqrt{3}, \sqrt{3}]$   the union of the three Hill intervals
 $[-\sqrt{3}, -1], [-1, 1], [1, \sqrt{3}]$. Moreover  $F(x) < -1$ for $|x| > \sqrt{3}$ and $F(x) \to - \infty$ as $x \to \infty$.

 We establish inequality (a) of  proposition \ref{prop: period asymptotics}    in two steps. 
 
\begin{figure}%
    \centering
    {{\includegraphics[width=4cm]{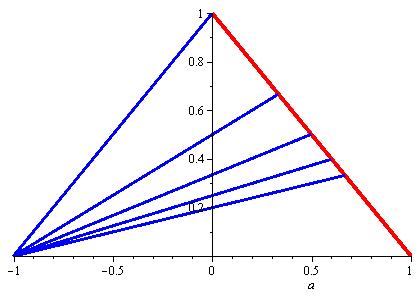} }}%
    \caption{The red line segment can be  parameterized by $(\mu,1-\mu)$ and represents $G$'s with Hill intervals $[-u, -1], [-1,1], [1, u]$. The blue line segments
    are parameterized by $\tau$ and represent those $G$'s whose  Hill interval  is
    fixed to be $[-u, u]$ with $u > \sqrt{3}$.}
    \label{Fol}
\end{figure}

 STEP 1. We show that   $Cost_y$ decreases monotonically along
 line segments in   $Diam_+$  connecting   $(-1,0)$ to points of $Leg 3$, 
 decreasing in the direction of $Leg 3$.  See   figure \ref{Fol} for this picture.

 STEP 2.  We show, with the help of numerics, that the absolute minimum of $Cost_y$ restricted to $Leg 3$
 occurs in the interior of the $Leg$ and is larger than $Cost_{0, y} (\infty)$. 
 
 \vskip .3cm
 
{\it The line segment family.} 
As in the  argument previously  used to get rid of points  near $(1,0)$ in $DIAM_-$,
these lines are determined by the condition that the far bound of integration $u$ 
(now characterized by $G(u) = -1$ since $b >0$) is constant along each line.
Parameterize $DIAM_+$ by coordinates $\tau, u$ with  $u > x_*$ and $-1 < \tau \le \frac{ |F(u)|-1}{|F(u)|+1 }$
according to 
\begin{eqnarray}
a & = \tau \\
b & = \frac{ \tau +1}{|F(u)|}.
\label{paramDiamond} 
\end{eqnarray}
Freezing $u$ and varying $\tau$ defines a line segment in the diamond which lies on the line whose slope is $1/|F(u)|$
passing through $(a,b) = (-1,0)$.  This line segment parameterizes    those $G = a + b F$'s  in $DIAM_+$  whose Hill interval is $[-u, u]$.
See the red line in figure \ref{Fol}.
Indeed since $G(x) = a + b F(x)$ we have that $G(u) = \tau - (\tau + 1) = -1$.
When $\tau = -1$ all these lines pass through $(-1,0)$.
They each reach $Leg 3$ ($a + b =1$) when  $\tau =  \frac{ |F(u)|-1}{|F(u)|+1 }$,
and their various endpoints sweep out $Leg 3$ as we vary $u$, with $u \to \infty$ 
being the limit of the   line segment  to $Leg 2$ ($b = 0$). 

  \vskip .3cm 
  
\subsubsection{Completing Step 1}

To complete step 1, we will prove that the cost  function $Cost(a, b)$ is a strictly  monotone decreasing function of $\tau$ relative to the $(\tau, u)$ parameterization
of $DIAM_+$ described by   equations (\ref{paramDiamond}).    In preparation for this computation observe that
   \begin{equation}\label{eq: costy}
   \begin{split}
Cost_y (a, b)  = &2 \int_0 ^u \frac{1 - G(x)}{ \sqrt{ 1- G(x)^2}} dx\\
 = & 2 \int_0 ^u \frac{\sqrt{1 - G(x)}}{ \sqrt{ 1+ G(x)}} dx,
\end{split}
   \end{equation}
and that, with points $(a,b) \in Diam_+$ parameterized by equation (\ref{paramDiamond})
we have
\begin{equation}
   \begin{split}
   G(x)  = & \tau + ( \frac{\tau + 1}{|F(u)|}) F(x) \\ 
  = & \frac{F(x)}{|F(u)|}  + \frac{\tau}{|F(u)|} (|F(u)| +  F(x)) ,
\end{split}
   \end{equation}
from which it follows that
\begin{eqnarray*}
\frac{d G}{d \tau} &= & \frac{1}{|F(u)|}(F(x)-F(u)),
  \end{eqnarray*}
  where we have used $F(u) < 0$.
  
 In the $(u, \tau)$ representation of points in $DIAM_+$, the $\tau$-lines 
 are lines of constant $u$ so we can   differentiate $Cost_y$ with respect to $\tau$
  by   differentiating the right hand side of equation (\ref{eq: costy}) under the integral sign.
    We get
    $$\frac{d}{d \tau} Cost_y( a(\tau, u), b(\tau, u)) = 2\int_0 ^u  \frac{d}{d \tau}  \frac{\sqrt{1 - G(x)}}{ \sqrt{ 1+ G(x)}}  dx,$$
   where $G = a + b F(x)$ and  $a = a(\tau, u), b = b(\tau, u)$ are given by
   equations (\ref{paramDiamond}). 
    Now then, the derivative of the integrand with respect to $\tau$ is  
   \begin{eqnarray*}
   \frac{d}{d \tau}  \frac{\sqrt{1 - G(x)}}{ \sqrt{ 1+ G(x)}}  & = &  \frac{-1}{ (1 - G(x))^{\frac{1}{2}}(1 + G(x))^{\frac{3}{2}}} \frac{dG}{d \tau} \\ 
   & = & \frac{1}{|F(u)|}\frac{F(u) - F(x)}{ (1 - G(x))^{\frac{1}{2}}(1 + G(x))^{\frac{3}{2}}} <0,  \\
   \end{eqnarray*}
the last inequality, the `$ < 0$',  holds  on the open interval $(0, u)$
as a result of the fact that  $F(u) < F(x)$ for $0 < x < u$. 
 It follows that the derivative of $Cost_y$ with respect to $\tau$ is strictly negative.

\vskip .3cm

\subsubsection{Completing Step 2}  

  $Cost_y$ extends  continuously to $Leg 3$, the  red line segment in figure \ref{Fol}. 
This extension is given by the same integral representation (\ref{eq: costy}) for $Cost_y$, with the upper bound of integration
 $u(a,b) $ being  the   positive solution to $a + b F(x) = -1$.  
What makes these $G$'s corresponding to    boundary points   $a + b =1$ on $Leg 3$ 
qualitatively  different than those $G$'s for $(a,b)$ in the   interior of the diamond is that the    $G$'s  for points on $Leg 3$ satisfy 
$G(1) = 1$ so   have three Hill regions:    $[-1,1]$,   $[-u, 1]$ and $[1, u]$ with $[-1,1]$ heteroclinic and the other two homoclinic.   
Parameterize $Leg 3$  by $a = \mu,  b = 1 - \mu$, for $0 < \mu < 1$.
Denote the continuous extension of $Cost_y$ to $Leg 3$ as `$Cost_{bdry}$'.
It forms a continuous function of   $\mu$, $1 \le \mu < 0$ upon setting  $a = \mu$ and $b = 1- \mu$ and  is given
  by  
$$Cost_{bdry} (\mu): =  2 \int_0 ^{u(\mu, 1- \mu)}\frac{\sqrt{1-G_{\mu}(x)} }{\sqrt{1+G_{\mu}(x)}}dx,$$
\beq 
G_{\mu} (x) :=   \mu + (1-\mu) F(x). 
\label{eq: Gmu} 
\eeq
Because of the nature of the Hill intervals of these $G_{\mu}$'s the  integral for
$Cost_{bdry} (\mu)$ incorporates the contributions of {\bf two} critical geodesics for $G_{\mu}$,
one being heteroclinic with interval $[-1,1]$ and the other being homoclinic with interval $1 \le x \le u: = u(\mu, 1-\mu)$.
Consequently we can write
$$Cost_{bdry} (\mu): = Cost_{heter} (\mu)  + Cost_{homoc} (\mu) , 0 < \mu < 1,$$
where the first integral goes from $0$ to $1$ and the second from $1$ to $u(\mu, 1- \mu)$. 
If we take the limit $\mu \to 0$ corresponding to $G_0(x) = F(x)$ , then $u \to \sqrt{3}$ from above
and we find that 
$$\lim_{\mu \to 0} Cost_{bdry} (\mu) = Cost_{0, y} (\infty) + Cost_{homoc}(\infty). $$
The  last term $Cost_{homoc}(\infty)$ is the integral from    $1$ to $\sqrt{3}$ of the integrand of equation (\ref{eq: costy} )
except with $G$ replaced by $F$. Expressed this way this last term is   clearly positive.     {\bf  It is this last positive jump of} $Cost_{homoc}(\infty)$ {\bf which  gives us our big advantage along $Leg 3$, making the completion of the proof possible.} 

\vskip .3cm

To finish the proof, we need to prove the following inequalities:

a) $Cost_{0, y} (\infty) < 0.58$.

\vskip .3cm

b) $0.58 < Cost_{bdry} (\mu)$ for all $\mu$ with $0 < \mu < 1$.

\vskip .3cm

To prove these inequalities we introduce auxilary functions
\begin{equation*}
\begin{split} 
F_{0}(x) &:= 1-\frac{(x^2-1)^2}{6}(1+2(\frac{x}{\sqrt{3}})^{14}), \\
 F_{1}(x) &:= 1-\frac{(x^2-1)^2}{6}(2+(\frac{x}{\sqrt{3}})^{14}). \\
\end{split}
\end{equation*}
$F_0(x)$ corresponds to $P(x) = (\frac{x}{\sqrt{3}})^{14}$ while $F_1(x)$ corresponds to  $P(x) = 1$
in   definition \ref{def: specific class} specifying the class of polynomials, the
definition  immediately following the statement of  theorem \ref{main}.

a)  Since $(\frac{x}{\sqrt{3}})^{14} \leq  P(x) \leq 1$ if $| x| < \sqrt{3}$ according to definition \ref{def: specific class} we have that
\begin{equation}\label{first-ineq}
F_{1}(x) \le F(x) \le F_{0} (x), \text{ for } | x| \le \sqrt{3}.
\end{equation}
The function $f \mapsto \sqrt{\frac{1-f}{1 + f}}$ is strictly monotone decreasing on the interval $[-\sqrt{3},\sqrt{3}]$
 from which it follows that 
 $$\sqrt{\frac{1-F_{0}(x)}{1+F_{0}(x)}} \le \sqrt{\frac{1-F(x)}{1+F(x)}} \le  \sqrt{\frac{1-F_{1}(x)}{1+F_{1}(x)}}, x \in [\sqrt{3}, \sqrt{3}].$$
 Using the upper bound for $x \in [-1,1]$ we get
 $$Cost_{0, y} (\infty) =  2  \int_{0}^1\sqrt{\frac{1-F(x)}{1+F(x)}}dx \le  2\int_{0}^1\sqrt{\frac{1-F_{1}(x)}{1+F_{1}(x)}}dx. $$
A numerical integration yields  
$$  2\int_0^1\sqrt{\frac{1-F_{1}(x)}{1+F_{1}(x)}}dx = 0.5790109314  , $$
establishing (a).

b)  Since $1 \leq P(x) \leq (\frac{x}{\sqrt{3}})^{14}$ if $ \sqrt{3} < | x| $ according to the specifications
of definition  \ref{def: specific class}  we have  
\begin{equation}\label{second-inq}
F_0(x) \leq F(x) \le F_{1} (x), \text{ for } \sqrt{3} \le |x| \le  \infty.
\end{equation}
We consider the family $G_{\mu}$ 
as in equation (\ref{eq: Gmu}).  Its Hill region has the form  $[-u(\mu), u(\mu)]$ where $\sqrt{3} \leq u(\mu)$,
$\mu \to u(\mu)$ is an invertible function.  (As before   the slope of the line joining $(-1,0)$ 
to $(\mu, 1-\mu)$ on $Leg 3$ is proportional to $1/|F(u)|$.)  We find $\mu(u)$ by solving   $G_{\mu(u)} (u) = -1$.
We get  $\mu(u) = \frac{-1-F(u)}{1-F(u)}$. Plug  $\mu(u)$ in $G_{\mu(u)}(x)$ and after a simplification we find
$$ Cost_{bdry}(u) = 2\int_0^u \sqrt{\frac{1-F(x)}{F(x)-F(u)}}dx.$$

Set 
$$ Lowbound(u) = 2\int_0^{\sqrt{3}} \sqrt{\frac{1-F_0(x)}{F_0(x)-F_0(u)}}dx + 2\int_{\sqrt{3}}^u \sqrt{\frac{1-F_1(x)}{F_0(x)-F_0(u)}}dx.$$
We claim that
$$Lowbound(u) <  Cost_{bdry}(u).$$
Indeed $1-F_0(x) < 1- F(x)$ on $[0,\sqrt{3}]$ and
$1 - F_1 (x) < 1- F(x)$ on $[\sqrt{3},u]$ by  inequality   (\ref{second-inq}).
And 
$ {\frac{1}{F_0(x)-F_0(u)}} < \frac{1 }{F(x)-F(u)}$ on the entire interval $[0, u]$.
This last inequality holds because    
\begin{equation*}
\begin{split}
F(x)-F(u) & \leq  F_0(x)-F_0(u), \text{ for } 0 \le x \le  u. \\
\end{split}
\end{equation*}  
which in turn  is true since it  is equivalent to 
\begin{equation*}
\begin{split}
  (x^2-1)^2((\frac{x}{\sqrt{3}})^{14} - P(x) ) & \leq (u^2-1)^2((\frac{u}{\sqrt{3}})^{14} - P(u) ), \text{ for } 0 \le x \le  u,
\end{split}
\end{equation*}
which is seen to hold upon using the properties of \ref{def: specific class}.
Indeed, according to these properties 
 $1 \leq P(u) \leq (\frac{u}{\sqrt{3}})^{14}$  since $u >\sqrt{3}$, so  the right side of the inequality is positive. 
 For $\sqrt{3} < x < u$ the left hand side is monotone increasing, reaching a maximum when $x =u$ and
 the inequality becomes equality. For $0 < x <\sqrt{3}$ the left  hand side  of the inequality
 is negative, again by one of the properties of \ref{def: specific class}.

\begin{figure}%
    \centering
    {{\includegraphics[width=8cm]{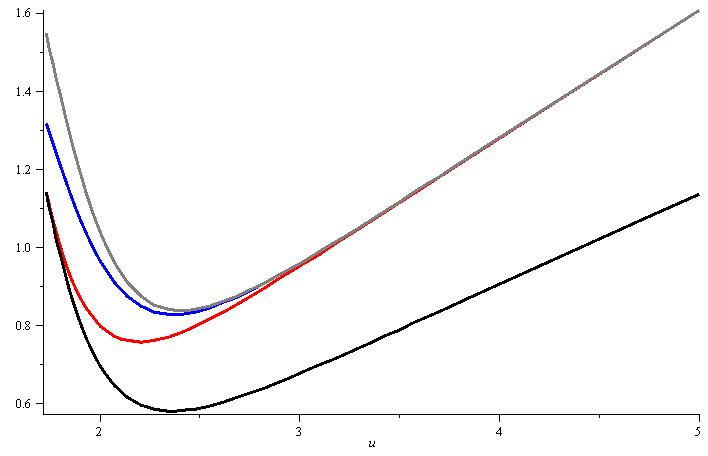} }}%
    \caption{On the panel graph of $Lowbound(u)$ being a lower bound for the graph of the $Cost_{bdry}(u)$ for $F_0(x)$,  $F_1(x)$ and $F(x)$ 
with $W(x)= \frac{1}{6}(1+ (\frac{x}{\sqrt{3}})^2 + (\frac{x}{\sqrt{3}})^{14})$.}
    \label{Cost-f-image}
\end{figure}

\begin{figure}%
    \centering
    {{\includegraphics[width=8cm]{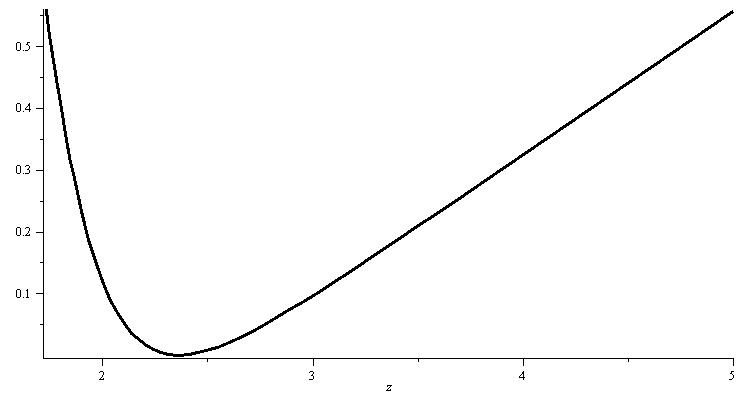} }}%
    \caption{On the right panel the graph of $Lowbound(u)$, where we can see that it is convex.}
    \label{Cost-f-image-1}
\end{figure}

Using Maple's built-in numerical integrator, we performed these integrals 
defining $Lowbound(u)$ to a tolerance of $10^{-10}$ and plotted the results.
See figure \ref{Cost-f-image} for a plot of the function $Lowbound(u) -.58$ versus $u$ for $u \ge \sqrt{3}$.  The plot shows
the function is positive and convex yielding that   $0.58 < Lowbound(u)$.   
(One can verify directly that   $Lowbound(u) \to \infty$ as $u \to \infty$.)   See also figure \ref{Cost-f-image-1}.

Now we are ready to finish the proof of proposition  \ref{prop: period asymptotics}.
  
\subsubsection{Proof of proposition  \ref{prop: period asymptotics}}

\begin{proof}
Inequality (a) of proposition \ref{prop: period asymptotics}  follows immediately from the lemma.
To repeat the argument:  Move from any interior $(a,b)$ point
to the boundary by moving along a $\tau$-line until you hit a point on  $Leg 3$.  
Since $Cost_y$ decreases monotonically along the  $\tau$-lines the value limited to on $Leg 3$,
whatever it be, is less than its  original value. Once  on   $Leg 3$,   $Cost_y$ is greater than
$Cost_0 (\infty)$.   
\end{proof}

\vskip .3cm

\subsection{Proof theorem \ref{main}}
 \label{subsec: end proof main} 
\begin{proof}
We show how   proposition \ref{prop: period asymptotics} completes the proof of Theorem \ref{main}.  
   All points of $DIAM_+$ are excluded by (a) of the proposition and the following logic.
   If a competing geodesic coming from $DIAM_+$ shares endpoints with $c_0$ for
   some $T$ then it shares y-endpoint values:  $\Delta y_0 (T) = \Delta y(a,b)$.
   By (a), $Cost_y (a,b) > Cost_{0, y} (\infty)$.  By lemma \ref{lem: cost at infinity}  $Cost_{0,y}(\infty) > Cost_{0,y}(T)$,
   so that  $\Delta t(a,b) - \Delta y(a,b) > T - \Delta y_0 (T)$. Thus $\Delta t (a,b) > T$:
   the competing geodesic is longer.
   
   All competing geodesics coming from points of $DIAM_-$ are excluded using (b) of the proposition and a similar logic.
   Since the endpoint conditions (equations(\ref{endpoints})) must hold for all sufficiently large $T$
   and since $\Delta z_0 (T), \Delta y_0 (T) \to \infty$ with $T$ we can restrict ourselves to
   $(a,b)$ in an arbitrarily small neighborhood of $Leg 1$ or $Leg 2$, since only here in $DIAM_-$ do the (a,b)
   periods blow up.  But then $Cost_z (a,b) \to \infty$ as we approach either leg.  In particular,
   for $(a,b)$ close to any point on either Leg we have $Cost_z (a,b) >> Cost_{0,z}(\infty) > Cost_{0,z}(T)$.
   Now the same logic and  inequalities as in the last two sentences
   of the previous paragraph carry through with $y$ replaced by $z$.  
\end{proof}

\begin{appendix} 

\section{Lagrangian approach to geodesic equations.}\label{lag-appro}

\begin{figure}
    \centering
    {{\includegraphics[width=10cm]{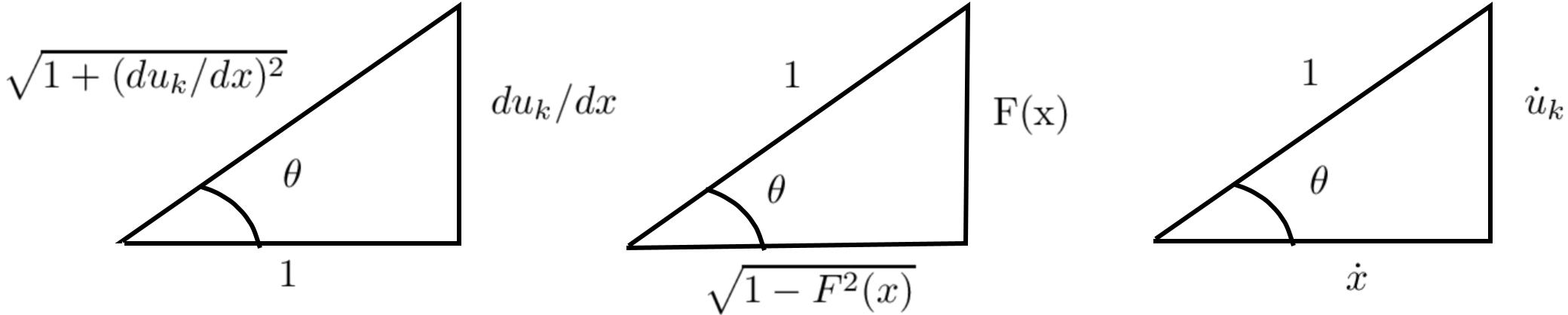} }}%
    \caption{The relation between the variables and the angle the curve makes with our orthonormal frame.}     \label{triangles}
\end{figure}

We derive the geodesic equations on the jet space $J^k$ using a Lagrangian approach. 
We can interpret  the coordinate $u_k$   as $\frac{d^k u_0}{d x^k}$.  We  can do so along any
arc of any horizontal curve $\gamma(t)$ in $J^k$ which is nowhere tangent
to the vertical  $X_2 = \dd{}{u_k}$ direction.  For then
$\dot x \ne 0$ so that we can  take $x$ rather that $t$ as the independent variable 
parameterizing the curve
Any such curve is then the $k$-jet of the function $f(x) = u_0 (x)$
so that along the curve   $u_i = d^i u_0 /dx ^i$.
We can  rewrite the arc-length of such a horizontal curve on $J^k$ as follow
\begin{eqnarray*}
 \int \sqrt{\dot{x}^2+\dot{u}_k^2}dt & = & \int (\sqrt{1+(\frac{d^{k+1} u_0}{d x^{k+1}})^2})\frac{dx}{dt} dt \\ 
& := &  \int \sqrt{1+(\frac{d^{k+1} u_0}{d x^{k+1}})^2} dx.
\end{eqnarray*}
The last integrand is a Lagrangian depending on higher derivatives,
but is independent of $x$.  Using the Euler-Lagrange equation for higher order derivatives (see \cite{Gelfand-Fomin}, pages 40-42) we have
$$ \frac{d^{k+1}}{dx^{k+1}}(\frac{\frac{d^{k+1} u_0}{d x^{k+1}}}{\sqrt{1+(\frac{d^{k+1} u_0}{d x^{k+1}})^2}})= 0 \;\; \Longrightarrow \;\; \frac{\frac{d^{k+1} u_0}{d x^{k+1}}}{\sqrt{1+(\frac{d^{k+1} u_0}{d x^{k+1}})^2}} = a_0+a_1x+\dots +a_k x^k. $$ 
Define $F(x):= a_0+a_1x+\dots +a_k x^k$ to be this guaranteed polynomial.
Observe that the function $\beta \mapsto f =  \beta/\sqrt{1 + \beta^2}$
is invertible with inverse $f \mapsto f / \sqrt{1 - f^2}$ and set $f = F(x)$, 
$\beta = \frac{d^{k+1} u_0}{d x^{k+1}} = du_k/dx$
to obtain  
$$ \frac{d u_k}{d x} = \frac{d^{k+1} u_0}{d x^{k+1}} = \frac{F(x)}{\sqrt{1-F^2(x)}}.$$
This equation is the same as the one from proposition \ref{y-z as x function}. 
See figure \ref{triangles} for the relation between the various  dependent variables and 
the angle $\theta$ made by our geodesic relative to the frame $X_1, X_2$,
i.e. to $\cos(\theta) X_1 + \sin(\theta) X_2 = \dot \gamma$.
We see that  
 $\frac{dx}{dt} = \sqrt{1-F^2(x)}$ and   $\dot{u}_k=F(x)$. Using the equations from the Pfaffian system, 
 equation (\ref{Pfaff}), (that is, what it means to locally be a k-jet)  we find the rest of the geodesics equations (\ref{eq:horiz}).  

\eject

\section{Multiplication on $J^k$.}

$J^k$ forms  a $k+2$-dimensional Lie group with multiplication 
\begin{equation*}
 (x,u_k,\dots,u_0) \cdot (x',u'_k,\dots,u'_0)  = (x+x',u_k+u'_k,u_{k-1} + u'_{k-1} + u_kx',\dots,u_0+u'_0 + u_1x').
\end{equation*}
 The neutral element is $(0,0,\dots,0)$ and the inverse is 
\begin{equation*}
(x,u_k,\dots,u_0)^{-1} = (-x,-u_k,-u_{k-1} + x'u_k,\dots,-u_0 + xu_1-x^2u_2+\dots+(-1)^{k-1}u_kx^k).
\end{equation*}
One computes  s t that the   standard  frame   $\{ X_1,\dots,X_{k+2}\}$  for $J^k$ generated by $X_1, X_2$ (equation (\ref{frame})) 
satisfies 
\begin{equation*}
\begin{split}
X_1(g) =(L_{g})_*(e_1) ,\; \cdots ,\; X_{k+2}(g) = (L_{g})_*(e_{k+2}),\\
\end{split}
\end{equation*}
where $g = (x,u_k,\dots,u_0)$,   $(L_{g})_*$ is the push forward by  left translation by $g$, and $e_i$ is the canonical base on $\R^{k+2}$.
This computation exhibits the  $X_i$ as a basis for the   left-invariant vector fields   on $J^k$.

\section{Inverting   polynomials near infinity}

In getting rid of points $(a,b)$ near Leg 2  having $b<0$
 we used an estimate for the solution to a real polynomial equation
$p(x) = w$ for $w= 1/\eps$ very large.   See equation (\ref{eq: root asymptotics}).  Here we derive that estimate.

Write
\begin{eqnarray*} w  & = &p(x) \\
                                 & = & a_0 x^n + a_1 x^{n-1} + \ldots  + a_n ,
 \end{eqnarray*}
Consider $w>>1$ and suppose $a_0 > 0$. 
 Write the reciprocal of both sides to get 
 \begin{eqnarray*} \frac{1}{w}   & = &\frac{1}{a_0 x^n + a_1 x^{n-1} + \ldots  + a_n}  \\
                                 & = & \frac{1}{a_0 x^n( 1+  b_1 x^{-1} + \ldots  + b_n x^{-n})}, 
 \end{eqnarray*}                            
where $b_i =a_i/a_0$.  Rewrite this relation in terms
of the coordinates  $u =1/x,  v = 1/ w$
 about  infinity in the domain and range:
 $$v = \frac{1}{a_0} u^n [ 1 + b_1 u + \ldots b_n u^n]^{-1},$$
 which we can partially invert 
 $$ (a_0 v)^{1/n} = u [1 + b_1 u + \ldots b_n u^n ]^{-1/n}.$$
 The right hand side function of $u$ is a near-identity analytic transformation near $u=0$
 so is invertible.  Write $f(u)$ for its inverse:  $f(u) = u - \frac{1}{n} (b_1 u )^2 + \ldots$
 by the binomial theorem.  Then  $u = f((a_0 v)^{1/n})$. 
 Since $x = 1/u$ and $w =1/v$  this gives us that
  $x = 1/f((a_0/w)^{1/n}) = 1/[(a_0 /w )^{1/n} -  \frac{1}{n} (b_1 (a_0 /w )^{1/n} )^2 + \ldots ]$.
  Expanding out we obtain the desired: 
  $$x = \frac{w^{1/n}}{a_0 ^{1/n}} + \frac{1}{n} b_1 + O(\frac{1}{w}).$$
  
%
%
%
%

\end{appendix}

\end{document}